\documentclass[10pt,a4paper]{article}

\usepackage{amsmath}
\usepackage[psamsfonts]{amssymb}
\usepackage[psamsfonts]{eucal}
\usepackage{amsthm}

\begin{document}

\newcommand{\Author}{Charles-Michel Marle}
\newcommand{\Address}{Universit\'{e} Pierre et Marie Curie,
                       Institut de Math\'{e}matiques,\\
                       4, place Jussieu,
                       75252 Paris cedex 05, France}

\newcommand{\Title}{Differential calculus on a Lie algebroid
                    and~Poisson manifolds}

\newcommand{\field}[1]{\mathbb{#1}}
\newcommand{\CC}{\field{C}}
\newcommand{\KK}{\field{K}}
\newcommand{\NN}{\field{N}}
\newcommand{\RR}{\field{R}}
\newcommand{\ZZ}{\field{Z}}

\swapnumbers
\newtheorem{ndefi}[subsubsection]{Definition}
\newtheorem{ntheorem}[subsubsection]{Theorem}
\newtheorem{nlemma}[subsubsection]{Lemma}
\newtheorem{ncorol}[subsubsection]{Corollary}
\newtheorem{nprop}[subsubsection]{Proposition}
\newtheorem{nrmk}[subsubsection]{Remark}
\newtheorem{nrmks}[subsubsection]{Remarks}
\newtheorem{nex}[subsubsection]{Example}
\newtheorem{nexs}[subsubsection]{Examples}
\newtheorem{defis}[subsubsection]{Definitions}

\newenvironment{nproof}{\noindent{\bf Proof:\\}}{ \hfill\qed}

\def\id{\mathop{\rm id}\nolimits}
\def\ad{\mathop{\rm ad}\nolimits}

\begin{center}
\Large{\bf{\Title}}
\end{center}

\begin{center}
\large{\bf{\Author}}
\end{center}

\begin{center}
\Address
\end{center}

\vspace{1cm}

\begin{abstract}
A Lie algebroid over a manifold is a vector bundle over that
manifold whose properties are very similar to those of a tangent
bundle. Its dual bundle has properties very similar to those of a
cotangent bundle: in the graded algebra of sections of its
external powers, one can define an operator $d_E$ similar to the
exterior derivative. We present in this paper the theory of Lie
derivatives, Schouten-Nijenhuis brackets and exterior derivatives
in the general setting of a Lie algebroid, its dual bundle and
their exterior powers. All the results (which, for their most
part, are already known) are given with detailed proofs. In the
final sections, the results are applied to Poisson manifolds.
 \end{abstract}
\vskip 1cm
\centerline{To Jos\'e~Antonio~Pereira~da~Silva and Alan~Weinstein}
\vskip 1cm

\section{Introduction}
\label{Introduction} %GIVE SECTION A LABEL%
\setcounter{equation}{0} Lie groupoids and Lie algebroids were
introduced in symplectic geometry by A.~Weinstein \cite{Wein2,
Coste} and, independently, M.~Karasev \cite{Kara}, in view of the
symplectization of Poisson manifolds and applications to
quantization. They are now an active domain of research, with
applications in various parts of mathematics \cite{Wein3, CannW,
Lib1}. More specifically, Lie algebroids have applications in
Mechanics \cite{Lib2} and provide a very natural setting in which
one can develop the theory of differential operators such as the
exterior derivative of forms and the Lie derivative with respect
to a vector field. In such a setting, slightly more general than
that of the tangent and cotangent bundles to a smooth manifold and
their exterior powers, the theory of Lie derivatives extends, in a
very natural way, into the theory of the Schouten-Nijenhuis
bracket (first introduced in differential geometry by
J.A.~Schouten \cite{Schou} and developed by A.~Nijenhuis
\cite{Ni}). Other bidifferential operators such as the bracket of
exterior forms on a Poisson manifold, first discovered for Pfaff
forms by F.~Magri and C.~Morosi \cite{MaMo} and extended to forms
of all degrees by J.-L.~Koszul \cite{Ko2} appear in such a setting
as very natural: they are Schouten-Nijenhuis brackets for the Lie
algebroid structure of the cotangent bundle to a Poisson manifold.
\par\smallskip
We present in this paper the theory of Lie derivatives,
Schouten-Nijenhuis brackets and exterior derivatives in the
general setting of a Lie algebroid, its dual bundle and their
exterior powers. All the results (which, for their most part, are
already known) are given with detailed proofs. Most of these
proofs are the same as the classical ones (when the Lie algebroid
is the tangent bundle to a smooth manifold); a few ones are
slightly more complicated, because the algebra of sections of
exterior powers of the dual of a Lie algebroid is not locally
generated by its elements of degree $0$ and their differentials
(contrary to the algebra of exterior differential forms on a
manifold). In the final section, the results are applied to
Poisson manifolds.

\section{Lie algebroids}
\label{Lie algebroids} %GIVE SECTION A LABEL%
\setcounter{equation}{0} \ The concept of a Lie algebroid was
first introduced by J.~Pradines \cite{Pra}, in relation with Lie
groupoids.
\subsection{Definition and examples}
 A Lie algebroid over a manifold is a vector bundle
based on that manifold, whose properties are very similar to those
of the tangent bundle. Let us give its formal definition.

 \begin{ndefi}\sl
Let $M$ be a smooth manifold and $(E,\tau,M)$ be a vector bundle
with base $M$. A {\it Lie algebroid structure\/} on that bundle is
the structure defined by the following data:
 \begin{description}
 \item[--] a composition law $(s_1,s_2)\mapsto\{s_1,s_2\}$ on the
space $\Gamma(\tau)$ of smooth sections of that bundle, for which
$\Gamma(\tau)$ becomes a Lie algebra,

 \item[--] a smooth vector bundle map $\rho:A\to TM$, where $TM$
is the tangent bundle of $M$, which satisfies the following two
properties:
   \begin{description}
   \item{\rm (i)} the map $s\mapsto \rho\circ s$ is a Lie algebras
homomorphism from the Lie algebra $\Gamma(\tau)$ into the Lie
algebra $A^1(M)$ of smooth vector fields on $M$;

   \item{\rm(ii)} for every pair $(s_1,s_2)$
of smooth sections of $\tau$, and every smooth function
$f:M\to\RR$, we have the Leibniz-type formula,
 \begin{equation*}\{s_1,fs_2\}=f\{s_1,s_2\}+\bigl(i(\rho\circ
s_1)df\bigr)s_2\,.
 \end{equation*}
   \end{description}
 \end{description}

The vector bundle $(E,\tau,M)$ equipped with its Lie algebroid
structure will be called a {\it Lie algebroid\/} and denoted by
$(E,\tau,M,\rho)$; the composition law
$(s_1,s_2)\mapsto\{s_1,s_2\}$ will be called the {\it bracket\/}
and the map $\rho:E\to TM$ the {\it anchor\/} of the Lie algebroid
$(E,\tau,M,\rho)$.
\end{ndefi}

\begin{nrmks} \rm Let $(E,\tau,M,\rho)$ be a Lie algebroid.
\par\nobreak\smallskip\noindent
{\rm(i)}\quad{\it The Lie algebras homomorphisms
$s\mapsto\rho\circ s$ and $s\mapsto{\cal L}(\rho\circ s)$.}\quad
Since $s\mapsto\rho\circ s$ is a Lie algebras homomorphism, we
have, for every pair $(s_1,s_2)$ of smooth sections of $\tau$,
 \begin{equation*}[\rho\circ s_1,\rho\circ s_2]=\rho\circ\{s_1,s_2\}\,.
 \end{equation*}

For every smooth vector field $X$ on $M$, let ${\cal L}(X)$ be the
Lie derivative with respect to that vector field. We recall that
${\cal L}(X)$ is a derivation of $A=C^\infty(M,\RR)$, {\it i.e.},
that for every pair $(f,g)$ of smooth functions on $M$,
 \begin{equation*}{\cal L}(X)(fg)=\bigl({\cal L}(X)f\bigr)g+f\bigl({\cal
 L}(X)g\bigr)\,.
 \end{equation*}
We recall also that $X\mapsto{\cal L}(X)$ is a Lie algebras
homomorphism from the Lie algebra $A^1(M)$ of smooth vector fields
on $M$, into the Lie algebra
 $\mathop{\rm Der}\bigl(C^\infty(M,\RR)\bigr)$
of derivations of $C^\infty(M,\RR)$, equipped with the commutator
 \begin{equation*}(D_1,D_2)\mapsto[D_1,D_2]=D_1\circ
 D_2-D_2\circ D_1
 \end{equation*}
as composition law.

The map $s\mapsto {\cal L}(\rho\circ s)$, obtained by composition
of two Lie algebras homomorphisms, is a Lie algebras homomorphism,
from the Lie algebra $\Gamma(\tau)$ of smooth sections of the Lie
algebroid $(E,\tau,M,\rho)$, into the Lie algebra of derivations of
$C^\infty(M,\RR)$.
\par\smallskip\noindent
{\rm(ii)}\quad{\it Leibniz-type formulae\/.}\quad According to
Definition 2.1.1 we have, for any pair $(s_1,s_2)$ of smooth
sections of $\tau$ and any smooth function $f$ on $M$,
 \begin{equation*}\{s_1,fs_2\}=f\{s_1,s_2\}+\bigl(i(\rho\circ
 s_1)df\bigr)\,s_2\,.
 \end{equation*}
As an easy consequence of the definition, we also have
 \begin{equation*}\{fs_1,s_2\}=f\{s_1,s_2\}-\bigl(i(\rho\circ
 s_2)df\bigr)\,s_1\,.
 \end{equation*}
More generally, for any pair $(s_1,s_2)$ of smooth sections of
$\tau$ and any pair $(f_1,f_2)$ of smooth functions  on $M$, we
have
 \begin{equation*}\{f_1s_1,f_2s_2\}=f_1f_2\{s_1,s_2\}+f_1\bigl(i(\rho\circ
 s_1)df_2\bigr)s_2-f_2\bigl(i(\rho\circ
 s_2)df_1\bigr)s_1\,.
 \end{equation*}
Using the Lie derivative operators, that formula may also be
written as
 \begin{equation*}\{f_1s_1,f_2s_2\}=f_1f_2\{s_1,s_2\}
 +f_1\bigl({\cal L}(\rho\circ
 s_1)f_2\bigr)s_2-f_2\bigl({\cal L}(\rho\circ
 s_2)f_1\bigr)s_1\,.
 \end{equation*}
\end{nrmks}
\par\bigskip

Let us give some examples of Lie algebroids.

\begin{nexs}\rm\hfill
\par\smallskip\noindent
{\rm(i)}\quad{\it The tangent bundle\/.}\quad The tangent bundle
$(TM,\tau_M,M)$ of a smooth manifold $M$, equipped with the usual
bracket of vector fields as composition law and with the identity
map $\id_{TM}$ as anchor, is a Lie algebroid.
\par\smallskip\noindent
{\rm(ii)}\quad{\it An involutive distribution\/.}\quad Let $V$ be
a smooth distribution on a smooth manifold $M$, {\it i.e.}, a
smooth vector subbundle of the tangent bundle $TM$. We assume that
$V$ is involutive, {\it i.e.}, such that the space of its smooth
sections is stable under the bracket operation. The vector bundle
$(V,\tau_M|_V,M)$, with the usual bracket of vector fields as
composition law and with the canonical injection $i_V:V\to TM$ as
anchor, is a Lie algebroid. We have denoted by $\tau_M:TM\to M$
the canonical projection of the tangent bundle and by $\tau_M|_V$
its restriction to the subbundle $V$.
\par\smallskip\noindent
{\rm(iii)}\quad{\it A sheaf of Lie algebras\/.}\quad Let
$(E,\tau,M)$ be a vector bundle over the smooth manifold $M$ and
$(z_1,z_2)\mapsto[z_1,z_2]$ be a smooth, skew-symmetric bilinear
bundle map defined on the fibered product $E\times_M E$, with
values in $E$, such that for each $x\in M$, the fibre
$E_x=\tau^{-1}(x)$, equipped with the bracket $(z_1,z_2)\mapsto
[z_1,z_2]$, is a Lie algebra. We define the bracket of two smooth
sections $s_1$ and $s_2$ of $\tau$ as the section $\{s_1,s_2\}$
such that, for each $x\in M$,
$\{s_1,s_2\}(x)=\bigl[s_1(x),s_2(x)\bigr]$. For the anchor, we
take the zero vector bundle map from $E$ to $TM$. Then $(E,\tau,M)$
is a Lie algebroid of particular type, called a {\it sheaf of Lie
algebras\/} over the manifold $M$.

\par\smallskip\noindent
{\rm(iv)}\quad{\it A finite-dimensional Lie algebra\/.}\quad In
particular, a finite-dimensional Lie algebra can be considered as
a Lie algebroid over a base reduced to a single point, with the
zero map as anchor.

\par\smallskip\noindent
{\rm(v)}\quad{\it The Lie algebroid of a Lie groupoid\/.}\quad To
every Lie groupoid, there is an associated Lie algebroid, much
like to every Lie group there is an associated Lie algebra. It is
in that setting that Pradines \cite{Pra} introduced Lie algebroids
for the first time. For more informations about Lie groupoids and
their associated Lie algebroids, the reader is referred to
\cite{Mack, Coste, DazSon, Albert}.
\end{nexs}

\subsection{Locality of the bracket}

We will prove that the value, at any point $x\in M$, of the
bracket of two smooth sections $s_1$ and $s_2$ of the Lie
algebroid $(E,\tau,M,\rho)$, depends only on the jets of order $1$
of $s_1$ and $s_2$ at $x$. We will need the following lemma.

\begin{nlemma}\sl
 Let $(E,\tau,M,\rho)$ be a Lie algebroid, $s_1:M\to E$ a smooth
 section of $\tau$, and $U$ an open subset of $M$ on which $s_1$
 vanishes. Then for any other smooth section $s_2$ of
 $\tau$, $\{s_1,s_2\}$ vanishes on $U$.

\end{nlemma}

\begin{nproof} Let $x$ be a point in $U$. There exists a smooth
function $f:M\to \RR$, whose support is contained in $U$ such that
$f(x)=1$. The section $fs_1$ vanishes identically, since $s_1$
vanishes on $U$ while $f$ vanishes outside of $U$. Therefore, for
any other smooth section $s_2$ of $\tau$,
 \begin{equation*}0=\{fs_1,s_2\}=-\{s_2,fs_1\}=-f\{s_2,s_1\}-\bigl(i(\rho\circ
 s_2)df\bigr)s_1\,.\end{equation*}
So at $x$ we have
 \begin{equation*}f(x) \{s_1,s_2\}(x)=\bigl(i(\rho\circ
 s_2)df\bigr)(x)s_1(x)=0\,.\end{equation*}
Since $f(x)=1$, we obtain $\{s_1,s_2\}(x)=0$.
\end{nproof}

\begin{nprop}\sl
 Let $(E,\tau,M,\rho)$ be a Lie algebroid. The value $\{s,s'\}(x)$
 of the bracket of two smooth sections $s$ and $s'$ of $\tau$,
 at a point $x\in M$, depends only on the jets of order $1$ of
 $s$ and $s'$ at $x$. Moreover, if $s(x)=0$ and $s'(x)=0$, then
 $\{s,s'\}(x)=0$.

\end{nprop}

\begin{nproof} Let $U$ be an open neighbourhood of $x$ in $M$ on
which there exists a local basis $(\sigma_1,\ldots,\sigma_k)$ of
of smooth sectionsof  $\tau$. For any point $y\in U$,
$\bigl(\sigma_1(y),\ldots,\sigma_k(y)\bigr)$ is a basis of the
fibre $E_y=\tau^{-1}(y)$. Let $s_1$ and $s_2$ be two smooth
sections of $\tau$. On the open subset $U$, these two sections can
be expressed, in a unique way, as
 \begin{equation*}s_1=\sum_{i=1}^k f_i\sigma_i\,,\qquad
  s_2=\sum_{j=1}^k g_j\sigma_j\,,\end{equation*}
where the $f_i$ and $g_j$ are smooth functions on $U$.
\par
By Lemma 2.2.1, the values of $\{s_1,s_2\}$ in $U$ depend only on
the values of $s_1$ and $s_2$ in $U$. Therefore we have in $U$
 \begin{equation*}\{s_1,s_2\}=\sum_{i,j}\Bigl(f_ig_j\{\sigma_i,\sigma_j\}
 +f_i\bigl({\cal L}(\rho\circ \sigma_i)g_j\bigr)\sigma_j
 -g_j\bigl({\cal L}(\rho\circ
\sigma_j)f_i\bigr)\sigma_i\Bigr)\,.\end{equation*} This expression
proves that the value of $\{s_1,s_2\}$ at $x$ depends only on the
$f_i(x)$, $df_i(x)$, $g_j(x)$ and $dg_j(x)$, that means on the
jets of order $1$ of $s_1$ and $s_2$ at $x$.
\par
If $s_1(x)=0$, we have, for all $i\in\{1,\ldots,k\}$,  $f_i(x)=0$,
and similarly if $s_2(x)=0$, we have, for all
$j\in\{1,\ldots,k\}$, $g_j(x)=0$. The above expression shows then
that $\{s_1,s_2\}(x)=0$.
\end{nproof}
\

\section{Exterior powers of vector bundles}
\label{Exterior powers} %GIVE SECTION A LABEL%
\setcounter{equation}{0} \

We recall in this section some definitions and general properties
related to vector bundles, their dual bundles and exterior powers.
In a first subsection we recall some properties of graded
algebras, graded Lie algebras and their derivations. The second
subsection applies these properties to the graded algebra of
sections of the exterior powers of a vector bundle. For more
details the reader may look at the book by Greub, Halperin and
Vanstone \cite{GrHV}. The reader already familiar with this
material may skip this section, or just look briefly at the sign
conventions we are using.

\subsection{Graded vector spaces and graded algebras}

\begin{defis}\sl\hfill
\par\nobreak\smallskip\noindent
{\rm(i)}\quad An {\it algebra\/} is a vector space $A$ on the
field $\KK=\RR$ or $\CC$, endowed with a $\KK$-bilinear map called
the {\it composition law\/},
 \begin{equation*}A\times A\to A\,,\quad (x,y)\mapsto xy\,,\quad \hbox{where}\
 x\ \hbox{and}\ y\in A\,.\end{equation*}
\par\smallskip\noindent
{\rm(ii)}\quad An algebra $A$ is said to be {\it associative\/} if
its composition law is associative, {\it i.e.}, if for all $x$,
$y$ and $z\in A$,
 \begin{equation*}x(yz)=(xy)z\,.\end{equation*}
\par\smallskip\noindent
{\rm(iii)}\quad A vector space $E$ on the field $\KK=\RR$ or $\CC$
is said to be {\it $\ZZ$-graded\/} if one has chosen a family
$(E^p\,,\ p\in\ZZ)$ of vector subspaces of $E$, such that
 \begin{equation*}E=\oplus_{p\in\ZZ}E^p\,.\end{equation*}
For each $p\in\ZZ$, an element $x\in E$ is said to be {\it
homogeneous of degree $p$ \/} if $x\in E^p$, and the vector
subspace $E^p$ of $E$ is called the {\it subspace of elements
homogeneous of degree $p$\/}.
\par\smallskip\noindent
{\rm(iv)}\quad Let $E=\oplus_{p\in\ZZ}E^p$ and
$F=\oplus_{p\in\ZZ}F^p$ be two $\ZZ$-graded vector spaces on the
same field $\KK$. A $\KK$-linear map $f:E\to F$ is said to be {\it
homogeneous of degree $d$\/} (with $d\in\ZZ$) if for each
$p\in\ZZ$,
  \begin{equation*}f(E^p)\subset F^{p+d}\,.\end{equation*}
\par\smallskip\noindent
{\rm(v)}\quad An algebra $A$ is said to be $\ZZ$-graded if
$A=\oplus_{p\in\ZZ}A^p$ is $\ZZ$-graded as a vector space and if
in addition, for all $p$ and $q\in \ZZ$, $x\in A^p$ and $y\in
A^q$,
  \begin{equation*}xy\in A^{p+q}\,.\end{equation*}
\par\smallskip\noindent
{\rm(vi)}\quad A $\ZZ$-graded algebra $A=\oplus_{p\in\ZZ}A^p$ is
said to be {\it $\ZZ_2$-commutative\/} (resp., {\it
$\ZZ_2$-anticommutative\/}) if for all $p$ and $q\in\ZZ$, $x\in
A^p$ and $y\in A^q$,
 \begin{equation*}xy=(-1)^{pq}yx\,,\qquad\hbox{(resp.,}\qquad
 xy=-(-1)^{pq}yx\,\hbox{)}\,.\end{equation*}

\end{defis}

\subsubsection{Some~properties~and~examples}
\par\smallskip\nobreak\noindent
{\rm(i)}\quad{\it Composition of homogeneous linear maps\/.}\quad
We consider three $\ZZ$-graded vector spaces,
$E=\oplus_{p\in\ZZ}E^p$, $F=\oplus_{p\in\ZZ}F^p$ and
$G=\oplus_{p\in\ZZ}G^p$, on the same field $\KK$. Let $f:E\to F$
and $g:F\to G$ be two linear maps, the first one $f$ being
homogeneous of degree $d_1$ and the second one $g$ homogeneous of
degree $d_2$. Then $g\circ f:E\to G$ is homogeneous of degree
$d_1+d_2$.
\par\smallskip\nobreak\noindent
{\rm(ii)}\quad {\it The algebra of linear endomorphisms of a
vector space\/.}\quad Let $E$ be a vector space and ${\cal
L}(E,E)$ be the space of linear endomorphisms of $E$. We take as
composition law on that space the usual composition of maps,
 \begin{equation*}(f,g)\mapsto f\circ g\,,\quad\hbox{with}\quad f\circ
 g(x)=f\bigl(g(x)\bigr)\,,\ x\in E\,.\end{equation*}
With that composition law, ${\cal L}(E,E)$ is an associative
algebra.
\par\smallskip\nobreak\noindent
{\rm(iii)}\quad{\it The graded algebra of graded linear
endomorphisms\/.}\quad We assume now that $E=\oplus_{p\in\ZZ}E^p$
is a $\ZZ$-graded vector space. For each $d\in\ZZ$, let $A^d$ be
the vector subspace of ${\cal L}(E,E)$ whose elements are the
linear endomorphisms $f:E\to E$ which are homogeneous of degree
$d$, {\it i.e.}, such that for all $p\in\ZZ$, $f(E^p)\subset
E^{p+d}$. Let $A=\oplus_{d\in\ZZ}A^d$. By using
property~3.1.2~(i), we see that with the usual composition of maps
as composition law, $A$ is a $\ZZ$-graded associative algebra.
\par\bigskip
Let us use property~3.1.2~(i) with $E=F=G$, in the following
definition.

\begin{ndefi}\sl Let $E=\oplus_{p\in\ZZ}E^p$ be a $\ZZ$-graded
vector space, $f$ and $g\in {\cal L}(E,E)$ be two homogeneous
linear endomorphisms of $E$ of degrees $d_1$ and $d_2$,
respectively. The linear endomorphism $[f,g]$ of $E$ defined by
 \begin{equation*}[f,g]=f\circ g-(-1)^{d_1d_2}g\circ f\,,\end{equation*}
which, by 3.1.2~(i), is homogeneous of degree $d_1+d_2$, is called
the {\it graded bracket\/} of $f$ and $g$.
\end{ndefi}

\begin{ndefi}\sl Let $A=\oplus_{p\in\ZZ}$ be a $\ZZ$-graded
algebra. Let $\theta:A\to A$ be a linear endomorphism of the
graded vector space $A$. Let $d\in\ZZ$. The linear endomorphism
$\theta$ is said to be a {\it derivation of degree $d$\/} of the
graded algebra $A$ if
\par\smallskip
\begin{description}
\item{\rm(i)} as a linear endomorphism of a graded vector space,
$\theta$ is homogeneous of degree $d$,

\item{\rm(ii)} for all $p\in\ZZ$, $x\in A^p$ and $y\in
A$,
 \begin{equation*}\theta(xy)=
 \bigl(\theta(x)\bigr)y+(-1)^{dp}x\bigl(\theta(y)\bigr)\,.\end{equation*}
\end{description}
\end{ndefi}

\begin{nrmk}\rm More generally, as shown by  Koszul \cite{Ko1}, for an
algebra $A$ equip\-ped with an involutive automorphism, one can
define two types of remarkable linear endomorphisms of $A$, the
{\it derivations\/} and the {\it antiderivations\/}. When
$A=\oplus_{p\in\ZZ}A^p$ is a $\ZZ$-graded algebra, and when the
involutive automorphism used is that which maps each $x\in A^p$
onto $(-1)^px$, it turns out that all nonzero graded derivations
are of even degree, that all nonzero graded antiderivations are of
odd degree, and that both derivations and antiderivations can be
defined as done in Definition 3.1.4. For simplicity we have chosen
to call {\it derivations\/} both the derivations and
antiderivations.
\end{nrmk}

\subsubsection{Some~properties~of~derivations}
Let $A=\oplus_{p\in\ZZ}A^p$ be a $\ZZ$-graded algebra.
\par\smallskip\nobreak\noindent
{\rm(i)}\quad{\it A derivation of degree $0$\/.}\quad
 For every
$p\in\ZZ$ and $x\in A^p$, we set
 \begin{equation*}\mu(x)=px\,.\end{equation*}
The map $\mu$, defined for homogeneous elements of $A$, can be
extended in a unique way as a linear endomorphism of $A$, still
denoted by $\mu$. This endomorphism is a derivation of degree $0$
of $A$.
\par\smallskip\nobreak\noindent
{\rm(ii)}\quad{\it The graded bracket of two derivations\,.}\quad
Let $\theta_1:A\to A$ and $\theta_2:A\to A$ be two derivations of
$A$, of degreeq $d_1$ and $d_2$, respectively. Their graded
bracket (Definition 3.1.3)
 \begin{equation*}[\theta_1,\theta_2]=\theta_1\circ\theta_2-(-1)^{d_1d_2}\theta_
2\theta_1\,,\end{equation*} is a derivation of degree $d_1+d_2$.

\begin{ndefi}\sl
A {\it $\ZZ$-graded Lie algebra\/} is a $\ZZ$-graded algebra
$A=\oplus_{p\in\ZZ}A^p$ (in the sense of 3.1.1~(v)), whose
composition law, often denoted by $(x,y)\mapsto [x,y]$ and called
the {\it graded bracket\/}, satisfies the following two
properties:
\par\smallskip
\begin{description}
\item{\rm(i)} it is $\ZZ_2$-anticommutative in the sense of 3.1.1~(vi),
{\it i.e.}, for all $p$ and $q\in\ZZ$, $P\in A^p$ and $Q\in A^q$,
 \begin{equation*}[P,Q]=-(-1)^{pq}[Q,P]\,,\end{equation*}

\item{\rm(ii)} it satisfies the {\it $\ZZ$-graded Jacobi identity\/}, {\it
i.e.}, for $p$, $q$ and $r\in\ZZ$, $P\in A^p$, $Q\in A^q$ and
$R\in A^r$,
 \begin{equation*}(-1)^{pr}\bigl[P,[Q,R]\bigr]+(-1)^{qp}\bigl[Q,[R,P]\bigr]
 +(-1)^{rq}\bigl[R,[P,Q]\bigr]=0\,.\end{equation*}
\end{description}
\end{ndefi}

\subsubsection{Examples~and~remarks}
\par\smallskip\nobreak\noindent
{\rm(i)}\quad{\it Lie algebras and $\ZZ$-graded Lie
algebras\/.}\quad A $\ZZ$-graded Lie algebra
$A=\oplus_{p\in\ZZ}A^p$ is not a Lie algebra in the usual sense,
unless $A^p=\{0\}$ for all $p\neq 0$. However, its subspace $A^0$
of homogeneous elements of degree $0$ is a Lie algebra in that
usual sense: it is stable under the bracket operation and when
restricted to elements in $A^0$, the bracket is skew-symmetric and
satisfies the usual Jacobi identity.

\par\smallskip\nobreak\noindent
{\rm(ii)}\quad{\it The graded Lie algebra associated to a graded
associative algebra\/.}\quad Let $A=\oplus_{p\in\ZZ}A^p$ be a
$\ZZ$-graded associative algebra, whose composition law is denoted
by $(P,Q)\mapsto PQ$. We define another composition law, denoted
by $(P,Q)\mapsto[P,Q]$ and called the {\it graded commutator\/};
we first define it for homogeneous elements in $A$ by setting, for
all $p$ and $q\in\ZZ$, $P\in A^p$ and $Q\in A^q$,
 \begin{equation*}[P,Q]=PQ-(-1)^{pq}QP\,;\end{equation*}
then we extend the definition of that composition law to all pairs
of elements in $A$ by bilinearity. The reader will easily verify
that with this composition law, $A$ is a graded Lie algebra. When
$A^p=\{0\}$ for all $p\neq 0$, we recover the well known way in
which one can associate a Lie algebra to any associative algebra.

\par\smallskip\nobreak\noindent
{\rm(iii)}\quad{\it The graded Lie algebra of graded
endomorphisms\/.}\quad Let $E=\oplus_{p\in\ZZ}E^p$ be a graded
vector space. For each $p\in\ZZ$, let $A^p\subset{\cal L}(E,E)$ be
the space of linear endomorphisms of $E$ which are homogeneous of
degree $p$, and let $A=\oplus_{p\in \ZZ}A^p$. As we have seen
in~3.1.2~(iii), when equipped with the composition of applications
as composition law, $A$ is a $\ZZ$-graded associative algebra. Let
us define another composition law on $A$, called the {\it graded
commutator\/}; we first define it for homogeneous elements in $A$
by setting, for all $p$ and $q\in\ZZ$, $P\in A^p$ and $Q\in A^q$,
 \begin{equation*}[P,Q]=PQ-(-1)^{pq}QP\,;\end{equation*}
then we extend the definition of that composition law to all pairs
of elements in $A$ by bilinearity. By using 3.1.8~(ii), we see
that $A$, with this composition law, is a $\ZZ$-graded Lie
algebra.

\par\smallskip\nobreak\noindent
{\rm(iv)}\quad{\it Various interpretations of the graded Jacobi
identity\/.}\quad Let $A=\oplus_{p\in\ZZ}A^p$ be a $\ZZ$-graded
Lie algebra. The $\ZZ$-graded Jacobi identity indicated in
Definition 3.1.7 can be cast into other forms, which better
indicate its meaning. Let us set, for all $P$ and $Q\in A$,
 \begin{equation*}\ad_PQ=[P,Q]\,.\end{equation*}
For each $p\in\ZZ$ and $P\in A^p$, $\ad_P:A\to A$ is a graded
endomorphism of $A$, homogeneous of degree $p$. By taking into
account the $\ZZ_2$-anticommutativity of the bracket, the reader
will easily see that the graded Jacobi identity can be written
under the following two forms:
\par\smallskip
{\it First form\/.}\quad For all $p$, $q$ and $r\in\ZZ$, $P\in
A^p$, $Q\in A^q$ and $R\in A^r$,
 \begin{equation*}\ad_P\bigl([Q,R]\bigr)
 =[\ad_PQ,R]+(-1)^{pq}[Q,\ad_PR]\,.
 \end{equation*}
This equality means that for all $p\in \ZZ$ and $P\in A^p$, the
linear endomorphism $\ad_P:A\to A$ is  a derivation of degree $p$
of the graded Lie algebra $A$, in the sense of 3.1.4.
\par\smallskip
{\it Second form\/.}\quad For all $p$, $q$ and $r\in \ZZ$, $P\in
A^p$, $Q\in A^q$ and $R\in A^r$,
 \begin{equation*}\ad_{[P,Q]}R=\ad_P\circ\ad_QR-(-1)^{pq}\ad_Q\circ\ad_PR=[\ad_P
,\ad_Q]R\,.\end{equation*} This equality means that for all $p$
and $q\in\ZZ$, $p\in A^p$ and $Q\in A^q$, the endomorphism
$\ad_{[P,Q]}:A\to A$ is the graded bracket (in the sense of 3.1.3)
of the two endomorphisms $\ad_P:A\to A$ and $\ad_Q:A\to A$. In
other words, the map $P\mapsto \ad_P$ is a $\ZZ$-graded Lie
algebras homomorphism from the $\ZZ$-graded Lie algebra $A$ into
the $\ZZ$-graded Lie algebra of sums of linear homogeneous
endomorphisms of $A$, with the graded bracket as composition law
(example 3.1.8~(iii)).
\par\smallskip
When $A^p=\{0\}$ for all $p\neq 0$, we recover the well known
interpretations of the usual Jacobi identity.

\subsection{Exterior powers of a vector bundle and of its dual}

In what follows all the vector bundles will be assumed to be
locally trivial and of finite rank; therefore we will write simply
{\it vector bundle\/} for {\it locally trivial vector bundle\/}.

\subsubsection{The~dual~of~a~vector~bundle} Let $(E,\tau,M)$ be a
vector bundle on the field $\KK=\RR$ or $\CC$. We will denote its
{\it dual bundle\/} by $(E^*,\pi,M)$. Let us recall that it is
a vector bundle over the same base manifold $M$, whose fibre
$E^*_x=\pi^{-1}(x)$ over each point $x\in M$ is the dual vector
space of the corresponding fibre $E_x=\tau^{-1}(x)$ of $(E,\tau,M)$,
{\it i.e.}, the space of linear forms on $E_x$ ({\it i.e.}, linear
functions defined on $E_x$ and taking their values in the field
$\KK$).
\par\smallskip
For each $x\in M$, the duality coupling $E^*_x\times E_x\to\KK$
will be denoted by
 \begin{equation*}(\eta,v)\mapsto\langle\eta,v\rangle\,.\end{equation*}

\subsubsection{The~exterior~powers~of~a~vector~bundle} Let
$(E,\tau,M)$ be a vector bundle of rank $k$. For each integer
$p>0$, we will denote by $(\bigwedge^pE,\tau,M)$ the $p$-th
external power of $(E,\tau,M)$. It is a vector bundle over $M$
whose fibre $\bigwedge^pE_x$, over each point $x\in M$, is the
$p$-th external power of the corresponding fibre $E_x=\tau^{-1}(x)$
of $(E,\tau,M)$. We recall that $\bigwedge^pE_x$ can be canonically
identified with the vector space of $p$-multilinear skew-symmetric
forms on the dual $E^*_x$ of $E_x$.
\par\smallskip
Similarly, for any integer $p>0$, we will denote by
$(\bigwedge^pE^*,\pi,M)$ the $p$-th external power of the
bundle $(E^*,\pi,M)$, dual of $(E,\tau,M)$.
\par\smallskip
For $p=1$, $(\bigwedge^1E,\tau,M)$ is simply the bundle
$(E,\tau,M)$, and similarly $(\bigwedge^1E^*,\pi,M)$ is simply
the bundle $(E^*,\pi,M)$. For $p$ strictly larger than the rank
$k$ of $(E,\tau,M)$, $(\bigwedge^pE,\tau,M)$ and
$(\bigwedge^pE^*,\pi,M)$ are the trivial bundle over $M$,
$(M\times\{0\},\tau_1,M)$, whose fibres are zero-dimensional
($\tau_1:M\times\{0\}\to M$ being the projection onto the first
factor).
\par\smallskip
For $p=0$, we will set
 $(\bigwedge^0E,\tau,M)=(\bigwedge^0E^*,\pi,M)=(M\times\KK,\tau_1,M)$,
where $\tau_1:M\times \KK\to M$ is the projection onto the first
factor.
\par\smallskip
Finally, we will consider that for $p<0$, $(\bigwedge^pE,\tau,M)$
and $(\bigwedge^pE^*,\pi,M)$ are the trivial bundle over $M$,
$(M\times\{0\},\tau_1,M)$. With these conventions,
$(\bigwedge^pE,\tau,M)$ and $(\bigwedge^pE^*,\pi,M)$ are defined
for all $p\in\ZZ$.
\par\medskip

\subsubsection{Operations in the graded vector spaces $\bigwedge
E_x$ and $\bigwedge E^*_x$}  Let $(E,\tau,M)$ be a vector bundle of
rank $k$, $(E^*,\pi,M)$ its dual and, for each $p\in\ZZ$,
$(\bigwedge^pE,\tau,M)$ and $(\bigwedge^pE^*,\pi,M)$ their
$p$-th external powers. We recall in this sections some operations
which can be made for each point $x\in M$, in the vector spaces
$\bigwedge^pE_x$ and $\bigwedge^pE^*_x$.
\par\smallskip
For each $x\in M$, let us consider the $\ZZ$-graded vector spaces
 \begin{equation*}\bigwedge
 E_x=\oplus_{p\in\ZZ}\bigwedge^pE_x\quad\hbox{and}\quad
 \bigwedge E^*_x=\oplus_{p\in\ZZ}\bigwedge^pE^*_x\,.\end{equation*}
We will say that elements in $\bigwedge E^*_x$ are (multilinear)
{\it forms\/} at $x$, and that elements in $\bigwedge E_x$ are
{\it multivectors\/} at $x$.
\par\smallskip\nobreak\noindent
{\rm(i)}\quad{\it The exterior product\/.}\quad Let us recall that
for each $x\in M$, $p$ and $q\in\ZZ$, $P\in\bigwedge^pE_x$ and
$Q\in\bigwedge^qE_x$, there exists $P\wedge
Q\in\bigwedge^{p+q}E_x$, called the {\it exterior product\/} of
$P$ and $Q$, defined by the following formulae.
\par\smallskip
\begin{description}
\item{--} If $p<0$, then $P=0$, therefore, for any
$Q\in\bigwedge^qE_x$, $P\wedge Q=0$. Similarly, if $q<0$, then
$Q=0$, therefore, for any $P\in\bigwedge^pE_x$, $P\wedge Q=0$.
\par\smallskip
\item{--} If $p=0$, then $P$ is a scalar ($P\in\KK$), and
therefore, for any $Q\in\bigwedge^qE_x$, $P\wedge Q=PQ$, the usual
product of $Q$ by the scalar $P$. Similarly, for $q=0$, then $Q$
is a scalar $(Q\in\KK)$, and therefore, for any
$P\in\bigwedge^pE_x$, $P\wedge Q=QP$, the usual product of $P$ by
the scalar $Q$.
\par\smallskip
\item{--} If $p\geq 1$ and $q\geq 1$, $P\wedge Q$, considered
as a $(p+q)$-multilinear form on $E^*_x$, is given by the formula,
where $\eta_1,\ldots,\eta_{p+q}\in E^*_x$,
\end{description}
 \begin{equation*}P\wedge Q(\eta_1,\ldots,\eta_{p+q})=\sum_{\sigma\in{\cal
 S}_{(p,q)}}\varepsilon(\sigma)
 P(\eta_{\sigma(1)},\ldots,\eta_{\sigma(p)})Q(\eta_{\sigma(p+1)},\ldots,
 \eta_{\sigma(p+q)})\,.\end{equation*}
We have denoted by ${\cal S}_{(p,q)}$ the set of permutations
$\sigma$ of $\{\,1,\ldots,p+q\,\}$ which satisfy
 \begin{equation*}\sigma(1)<\sigma(2)<\cdots<\sigma(p)\quad\hbox{and}\quad
 \sigma(p+1)<\sigma(p+2)<\cdots<\sigma(p+q)\,,\end{equation*}
and set
 \begin{equation*}\epsilon(\sigma)=\begin{cases}
                    1&\text{if $\sigma$ is even},\\
                    -1&\text{if $\sigma$ is odd}.
                    \end{cases}\end{equation*}
\par\smallskip
Similarly, let us recall that for each $x\in M$, $p$ and
$q\in\ZZ$, $\xi\in\bigwedge^pE^*_x$ and $\eta\in\bigwedge^qE^*_x$,
there exists $\xi\wedge \eta\in\bigwedge^{p+q}E^*_x$, called the
{\it exterior product\/} of $\xi$ and $\eta$. It is defined by the
formulae given above, the only change being the exchange of the
roles of $E_x$ and $E^*_x$.
\par\smallskip
The exterior product is associative and $\ZZ_2$-commutative: for
all $x\in M$, $p$, $q$ and $r\in\ZZ$, $P\in \bigwedge^pE_x$, $Q\in
\bigwedge^qE_x$ and $R\in \bigwedge^rE_x$,
 \begin{equation*}P\wedge(Q\wedge R)=(P\wedge Q)\wedge R\,,\qquad Q\wedge
 P=(-1)^{pq}P\wedge Q\,,\end{equation*}
 and similarly, for $\xi\in\bigwedge^pE^*_x$, $\eta\in
\bigwedge^qE^*_x$ and $\zeta\in \bigwedge^rE^*_x$,
 \begin{equation*}\xi\wedge(\eta\wedge \zeta)=(\xi\wedge \eta)\wedge
\zeta\,,\qquad \eta\wedge
 \xi=(-1)^{pq}\xi\wedge \eta\,.\end{equation*}
\par\medskip
For all $x\in M$, the exterior product extends, by bilinearity, as
a composition law in each of the graded vector spaces $\bigwedge
E_x$ and $\bigwedge E^*_x$. With these composition laws, these
vector spaces become $\ZZ$-graded associative and
$\ZZ_2$-commutative algebras.
\par\smallskip\nobreak\noindent
{\rm(ii)}\quad{\it The interior product of a form by a
vector\/.}\quad Let us recall that for each $x\in M$, $v\in E_x$,
$p\in\ZZ$, $\eta\in\bigwedge^pE^*_x$, there exists $i(v)\eta\in
\bigwedge^{p-1}E^*_x$, called the {\it interior product\/} of
$\eta$ by $v$, defined by the following formulae.
\par\smallskip
\begin{description}
\item{--} For $p\leq 0$, $i(v)\eta=0$, since
$\bigwedge^{p-1}E^*_x=\{0\}$.
\par\smallskip
\item{--} For $p=1$,
 \begin{equation*}i(v)\eta=\langle\eta,v\rangle\in\KK\,.\end{equation*}

\par\smallskip
\item{--} For $p>1$, $i(v)\eta$ is the $(p-1)$-multilinear form
on $E_x$ such that, for all $v_1,\ldots,v_{p-1}\in E_x$,
 \begin{equation*}i(v)\eta(v_1,\ldots,v_{p-1})=\eta(v,v_1,\ldots,v_{p-1})\,.\end
{equation*}
\end{description}
\par\medskip
For each $x\in M$ and $v\in E_x$, the map $\eta\mapsto i(v)\eta$
extends, by linearity, as a graded endomorphism of degree $-1$ of
the graded vector space $\bigwedge E^*_x$. Moreover, that
endomorphism is in fact a derivation of degree $-1$ of the
exterior algebra of $E^*_x$, {\it i.e.}, for all $p$ and
$q\in\ZZ$, $\zeta\in \bigwedge^pE^*_x$, $\eta\in\bigwedge^qE^*_x$,
 \begin{equation*}i(v)(\zeta\wedge\eta)
 =\bigl(i(v)\zeta\bigr)\wedge\eta+(-1)^p
 \zeta\wedge\bigl(i(v)\eta\bigr)\,.
 \end{equation*}
\par\smallskip\nobreak\noindent
{\rm(iii)}\quad{\it The pairing between $\bigwedge E_x$ and
$\bigwedge E^*_x$\/.}\quad Let $x\in M$, $p$ and $q\in \ZZ$,
$\eta\in\bigwedge^pE^*_x$ and $v\in\bigwedge^qE_x$. We set
 \begin{equation*}\langle\eta,v\rangle=\begin{cases}
                        0&\text{if $p\neq q$, or if $p<0$, or if
                        $q<0$},\\
                        \eta v&\text{if $p=q=0$}.
                        \end{cases}
 \end{equation*}
In order to define $\langle\eta,v\rangle$ when $p=q\geq 1$, let us
first assume that $\eta$ and $v$ are decomposable, {\it i.e.},
that they can be written as
 \begin{equation*}\eta=\eta_1\wedge\cdots\wedge\eta_p\,,\qquad
 v=v_1\wedge\cdots\wedge v_p\,,\end{equation*}
where $\eta_i\in E^*_x$, $v_j\in E_x$, $1\leq i,j\leq p$. Then we
set
 \begin{equation*}\langle\eta,v\rangle=\det\bigl(\langle\eta_i,v_j\rangle\bigr)\
,.\end{equation*} One may see that $\langle\eta,v\rangle$ depends
only on $\eta$ and $v$, not on the way in which they are expressed
as exterior products of elements of degree $1$. The map
$(\eta,v)\mapsto\langle\eta,v\rangle$ extends, in a unique way as
a bilinear map
 \begin{equation*}\bigwedge E^*_x\times\bigwedge E_x\to\KK,\quad\hbox{still
 denoted by}\quad
 (\eta,v)\mapsto\langle\eta,v\rangle\,,\end{equation*}
called the {\it pairing\/}. That map allows us to consider each
one of the two graded vector spaces $\bigwedge E^*_x$ and
$\bigwedge E_x$ as the dual of the other one.
\par\smallskip
Let $\eta\in\bigwedge^pE^*_x$ and  $v_1,\ldots,v_p$ be elements of
$E_x$. The pairing $\langle\eta,v_1\wedge\cdots\wedge v_p\rangle$
is related, in a very simple way, to the value of $\eta$,
considered as a $p$-multilinear form on $E_x$,  on the set
$(v_1,\ldots,v_p)$. We have
 \begin{equation*}\langle\eta,v_1\wedge\cdots\wedge
 v_p\rangle=\eta(v_1,\ldots,v_p)\,.\end{equation*}
\par\smallskip\nobreak\noindent
{\rm(iv)}\quad{\it The interior product of a form by a
multivector\/.}\quad For each $x\in M$ and $v\in E_x$,  we have
defined in 3.2.3~(ii) the interior product $i(v)$ as a derivation
of degree $-1$ of the exterior algebra $\bigwedge E^*_x$ of forms
at $x$. Let us now define, for each multivector $P\in\bigwedge
E_x$, the interior product $i(P)$. Let us first assume that $P$ is
homogeneous of degree $p$, {\it i.e.}, that $P\in\bigwedge^pE_x$.
\par\smallskip
\begin{description}
\item{--} For $p<0$, $\bigwedge^pE_x=\{0\}$, therefore $i(P)=0$.
\par\smallskip
\item{--} For $p=0$, $\bigwedge^0E_x=\KK$, therefore $P$ is a
scalar and we set, for all $\eta\in\bigwedge E^*_x$,
 \begin{equation*}i(P)\eta=P\eta\,.\end{equation*}
\par\smallskip
\item{--} For $p\geq 1$ and $P\in\bigwedge^pE_x$ decomposable,
{\it i.e.},
 \begin{equation*}P=P_1\wedge \cdots\wedge P_p\,,\quad \hbox{with}\quad
P_i\in E_x\,,
 \quad 1\leq i\leq p\,,\end{equation*}
we set
 \begin{equation*}i(P_1\wedge\cdots\wedge P_p)=i(P_1)\circ\cdots\circ
i(P_p)\,.\end{equation*} We see easily that $i(P)$ depends only of
$P$, not of the way in which it is expressed as an exterior
product of elements of degree $1$.
\par\smallskip
\item{--} We extend by linearity the definition of
$i(P)$ for all $P\in\bigwedge^pE_x$, and we see that $i(P)$ is a
graded endomorphism of degree $-p$ of the graded vector space
$\bigwedge E^*_x$. Observe that for $p\neq 1$, $i(P)$ is not in
general a derivation of the exterior algebra $\bigwedge E^*_x$.
\end{description}
\par\smallskip
Finally, we extend by linearity the definition of $i(P)$ to all
elements $P\in\bigwedge E_x$.
\par\smallskip\nobreak\noindent
{\rm(v)}\quad{\it The interior product by an exterior
product\/.}\quad It is easy to see that for all $P$ and
$Q\in\bigwedge E_x$,
 \begin{equation*}i(P\wedge Q)=i(P)\circ i(Q)\,.\end{equation*}
\par\smallskip\nobreak\noindent
{\rm(vi)}\quad{\it Interior product and pairing\/.}\quad For
$p\in\ZZ$, $\eta\in\bigwedge^{p}E^*_x$ and $P\in\bigwedge^pE_x$,
we have
 \begin{equation*}i(P)\eta=(-1)^{(p-1)p/2}\langle\eta,P\rangle\,.\end{equation*}
More generally, for $p$ and $q\in\ZZ$, $P\in\bigwedge^p(E_x)$,
$Q\in\bigwedge^q(E_x)$ and  $\eta\in\bigwedge^{p+q}(E^*_x)$,
 \begin{equation*}\bigl\langle
i(P)\eta,Q\bigr\rangle=(-1)^{(p-1)p/2}\langle\eta,
 P\wedge Q\rangle\,.\end{equation*}
This formula shows that the interior product by
$P\in\bigwedge^pE_x$ is $(-1)^{(p-1)p/2}$ times the transpose,
with respect to the pairing, of the exterior product by $P$ on the
left.

\subsubsection{The~exterior~algebra~of~sections} Let $(E,\tau,M)$ be
a vector bundle of rank $k$ on the field $\KK=\RR$ or $\CC$, over
a smooth manifold $M$, $(E^*,\pi,M)$ be its dual bundle and,
for each integer $p\geq 1$, let $(\bigwedge^pE,\tau,M)$ and
$(\bigwedge^pE^*,\pi,M)$ be their respective $p$-th exterior
powers.
\par\smallskip
For each $p\in\ZZ$, we will denote by $A^p(M,E)$ the space of
smooth sections of $(\bigwedge^pE,\tau,M)$, {\it i.e.}, the space
of smooth maps $Z:M\to\bigwedge^pE$ which satisfy
 \begin{equation*}\tau\circ Z=\id_M\,.\end{equation*}
Similarly, for each $p\in\ZZ$, we will denote by $\Omega^p(M,E)$
the space of smooth sections of $(\bigwedge^pE^*,\pi,M)$, {\it
i.e.}, the space of smooth maps $\eta:M\to\bigwedge^pE^*$ which
satisfy
 \begin{equation*}\pi\circ \eta=\id_M\,.\end{equation*}
Let us observe that $\Omega^p(M,E)=A^p(M,E^*)$.
\par
We will denote by $A(M,E)$ and $\Omega(M,E)$ the direct sums
 \begin{equation*}A(M,E)=\oplus_{p\in\ZZ}A^p(M,E)\,,\qquad
   \Omega(M,E)=\oplus_{p\in\ZZ}\Omega^p(M,E)\,.\end{equation*}
These direct sums, taken for all $p\in\ZZ$, are in fact taken for
all integers $p$ which satisfy $0\leq p\leq k$, where $k$ is the
rank of the vector bundle $(E,\tau,M)$, since we have
$A^p(M,E)=\Omega^p(M,E)=\{0\}$ for $p<0$ as well as for $p>k$.
\par\smallskip
For $p=0$, $A^0(M,E)$ and $\Omega^0(M,E)$ both coincide with the
space $C^\infty(M,\KK)$ of smooth functions defined on $M$ which
take their values in the field $\KK$.
\par\smallskip
Operations such as the exterior product, the interior product and
the pairing, defined for each point $x\in M$ in 3.2.3, can be
extended to elements in $A(M,E)$ and $\Omega(M,E)$.
\par\smallskip\nobreak\noindent
{\rm(i)}\quad{\it The exterior product of two sections\/.}\quad
For example, the exterior product of two sections $P$ and $Q\in
A(M,E)$ is the section
 \begin{equation*}x\in M\,,\quad x\mapsto (P\wedge Q)(x)=P(x)\wedge
Q(x)\,.\end{equation*} The exterior product of two sections $\eta$
and $\zeta\in\Omega(M,E)$ is similarly defined.
\par\smallskip
With the exterior product as composition law, $A(M,E)$ and
$\Omega(M,E)$ are $\ZZ$-graded associative and $\ZZ_2$-commutative
algebras, called the {\it algebra of multivectors\/} and the {\it
algebra of forms\/} associated to the vector bundle $(E,\tau,M)$.
Their subspaces $A^0(M,E)$ and $\Omega^0(M,E)$ of homogeneous
elements of degree $0$ both coincide with the usual algebra
$C^\infty(M,\KK)$ of smooth $\KK$-valued functions on $M$, with
the usual product of functions as composition law. We observe that
$A(M,E)$ and $\Omega(M,E)$ are $\ZZ$-graded modules over the ring
of functions $C^\infty(M,\KK)$.
\par\smallskip\nobreak\noindent
{\rm(ii)}\quad{\it The interior product by a section of
$A(M,E)$\/.}\quad For each $P\in A(M,E)$, the {\it interior
product\/} $i(P)$ is an endomorphism of the graded vector space
$\Omega(M,E)$. If $p\in\ZZ$ and $P\in A^p(M,E)$, the endomorphism
$i(P)$ is homogeneous of degree $-p$. For $p=1$, $i(P)$ is a
derivation of degree $-1$ of the algebra $\Omega(M,E)$.
\par\smallskip\nobreak\noindent
{\rm(iii)}\quad{\it The pairing between $A(M,E)$ and
$\Omega(M,E)$\/.}\quad The {\it pairing\/}
 \begin{equation*}(\eta,P)\mapsto\langle\eta,P\rangle\,,\quad \eta\in
 \Omega(M,E)\,,\quad P\in A(M,E)\,,\end{equation*}
is a $C^\infty(M,\KK)$-bilinear map, defined on $\Omega(M,E)\times
A(M,E)$, which takes its values in $C^\infty(M,\KK)$. \

\section{Exterior powers of a Lie algebroid and of its dual}
\label{Exterior powers 2} %GIVE SECTION A LABEL%
\setcounter{equation}{0} \ We consider now a Lie algebroid
$(E,\tau,M,\rho)$ over a smooth manifold $M$. We denote by
$(E^*,\pi,M)$ its dual vector bundle, and use all the notations
defined in Section~3. We will assume that the base field $\KK$ is
$\RR$, but most results remain valid for $\KK=\CC$. We will prove
that differential operators such as the Lie derivative and the
exterior derivative, which are well known for sections of the
exterior powers of a tangent bundle and of its dual, still exist
in this more general setting.

\subsection{Lie derivatives with respect to sections of a Lie
algebroid}

We prove in this subsection that for each smooth section $V$ of
the Lie algebroid $(E,\tau,M,\rho)$, there exists a derivation of
degree $0$ of the exterior algebra $\Omega(M,E)$, called the {\it
Lie derivative\/} with respect to $V$ and denoted by ${\cal
L}_\rho(V)$. When the Lie algebroid is the tangent bundle
$(TM,\tau_M,M,\id_{TM})$, we will recover the usual Lie derivative
of differential forms with respect to a vector field.

\begin{nprop}\sl Let $(E,\tau,M,\rho)$ be a Lie algebroid on a smooth
manifold $M$. For each smooth section $V\in A^1(M,E)$ of the
vector bundle $(E,\tau,M)$, there exists a unique graded
endomorphism of degree $0$ of the graded algebra of exterior forms
$\Omega(M,E)$, called the {\it Lie derivative\/} with respect to
$V$ and denoted by ${\cal L}_\rho(V)$, which satisfies the
following properties:
\begin{description}
\par\smallskip
\item{\rm(i)} For a smooth function
$f\in\Omega^0(M,E)=C^\infty(M,\RR)$,
 \begin{equation*}{\cal L}_\rho(V)f=i(\rho\circ V)df={\cal L}(\rho\circ
V)f\,,\end{equation*} where ${\cal L}(\rho\circ V)$ denotes the
usual Lie derivative with respect to the vector field $\rho\circ
V$;
\par\smallskip
\item{\rm (ii)} For a form $\eta\in\Omega^p(M,E)$ of degree $p>0$,
${\cal L}_\rho(V)\eta$ is the form defined by the formula, where
$V_1,\ldots,V_p$ are smooth sections of $(E,\tau,M)$,
 \begin{equation*}
 \begin{split}
 \bigl({\cal L}_\rho(V)\eta\bigr)(V_1,\ldots,V_p)
 &={\cal L}_\rho(V)\bigl(\eta(V_1,\ldots,V_p)\bigr)\\
 &\quad-\sum_{i=1}^p\eta(V_1,\ldots,V_{i-1},\{V,V_i\},V_{i+1},\ldots,V_p)\,.
 \end{split}
 \end{equation*}
\par\smallskip
\end{description}
\end{nprop}

\begin{nproof} Clearly (i) defines a function ${\cal L}_\rho(V)f\in
\Omega^0(M,E)=C^\infty(M,\RR)$. We see immediately that for $f$
and $g\in C^\infty(M,\RR)$,
 \begin{equation*}{\cal L}_\rho(V)(fg)=\bigl({\cal L}_\rho(V)f\bigr)g
 +f\bigl({\cal L}_\rho(V)g\bigr)\,.\eqno(*)\end{equation*}
Now (ii) defines a map $(V_1,\ldots V_p)\mapsto\bigl({\cal
L}_\rho(V)\eta\bigr)(V_1,\ldots,V_p)$ on $\bigl(A^1(M,E)\bigr)^p$,
with values in $C^\infty(M,\RR)$. In order to prove that this map
defines an element ${\cal L}_\rho(V)\eta$ in $\Omega^p(M,E)$, it
is enough to prove that it is skew-symmetric and
$C^\infty(M,\RR)$-linear in each argument. The skew-symmetry and
the $\RR$-linearity in each argument are easily verified. There
remains only to prove that for each function $f\in
C^\infty(M,\RR)$,
 \begin{equation*}\bigl({\cal L}_\rho(V)\eta\bigr)(fV_1,V_2,\ldots,V_p)=f
 \bigl({\cal
L}_\rho(V)\eta\bigr)(V_1,V_2,\ldots,V_p)\,.\eqno(**)\end{equation*}
We have
 \begin{equation*}
 \begin{split}
 \bigl({\cal L}_\rho(V)\eta\bigr)(fV_1,V_2,\ldots,V_p)
 &={\cal L}_\rho(V)\bigl(\eta(fV_1,V_2,\ldots,V_p)\bigr)\\
 &\quad-\eta(\{V,fV_1\},V_2,\ldots,V_p)\\
 &\quad-\sum_{i=2}^p\eta(fV_1,V_2,\ldots,V_{i-1},\{V,V_i\},V_{i+1},\ldots
 ,V_p)\,.
 \end{split}
 \end{equation*}
By using $(*)$, we may write
 \begin{equation*}
 \begin{split}
 {\cal L}_\rho(V)\bigl(\eta(fV_1,V_2,\ldots,V_p)\bigr)
 &={\cal L}_\rho(V)\bigl(f\eta(V_1,V_2,\ldots,V_p)\bigr)\\
 &=\bigl({\cal L}_\rho(V)f\bigr)\eta(V_1,V_2,\ldots,V_p)\\
 &\quad+f{\cal
  L}_\rho(V)\bigl(\eta(V_1,V_2,\ldots,V_p)\bigr)\,.
 \end{split}
  \end{equation*}
Using the property of the anchor, we also have
 \begin{equation*}\{V,fV_1\}=\bigl(i(\rho\circ
V)df\bigr)V_1+f\{V,V_1\}=\bigl({\cal
 L}_\rho(V)f\bigr)V_1+f\{V,V_1\}\,.\end{equation*}
Equality $(**)$ follows immediately.
\par
The endomorphism ${\cal L}_\rho(V)$, defined on the subspaces of
homogeneous forms, can then be extended, in a unique way, to
$\Omega(M,E)$, by imposing the $\RR$-linearity of the  map
$\eta\mapsto {\cal L}_\rho(V)\eta$.
\end{nproof}
\par\bigskip

Let us now introduce the $\Omega(M,E)$-valued exterior derivative
of a function. In the next section, that definition will be
extended to all elements in $\Omega(M,E)$.

\begin{ndefi}\sl Let $(E,\tau,M,\rho)$ be a Lie algebroid on a
smooth manifold $M$. For each function
$f\in\Omega^0(M,E)=C^\infty(M,\RR)$, we call {\it
$\Omega(M,E)$-valued exterior derivative\/} of $f$, and denote by
$d_\rho f$, the unique element in $\Omega^1(M,E)$ such that, for
each section $V\in A^1(M,E)$,
 \begin{equation*}\langle d_\rho f,V\rangle=\langle df,\rho\circ
V\rangle\,.\end{equation*}

\end{ndefi}

\begin{nrmk}\rm Let us observe that the transpose of the anchor
$\rho:E\to TM$ is a vector bundle map ${}^t\!\rho:T^*M\to E^*$. By
composition of that map with the usual differential of functions,
we obtain the $\Omega(M,E)$-valued exterior differential $d_\rho$.
We have indeed
 \begin{equation*}d_\rho f={}^t\!\rho\circ df\,.\end{equation*}
\end{nrmk}

\begin{nprop}\sl Under the assumptions of Proposition
4.1.1, the Lie derivative has the following properties.
\par\smallskip
{\rm 1.} For each $V\in A^1(M,E)$ and $f\in C^\infty(M,\RR)$,
 \begin{equation*}{\cal L}_\rho(V)(d_\rho f)=d_\rho\bigl({\cal
 L}_\rho(V)f\bigr)\,.\end{equation*}
\par\smallskip
{\rm 2.} For each $V$ and $W\in A^1(M,E)$, $\eta\in\Omega(M,E)$,
 \begin{equation*}i\bigl(\{V,W\}\bigr)\eta=\bigl({\cal L}_\rho(V)\circ i(W)-
 i(W)\circ{\cal L}_\rho(V)\bigr)\eta\,.\end{equation*}
\par\smallskip
{\rm 3.} For each $V\in A^1(M,E)$, ${\cal L}_\rho(V)$ is a
derivation of degree $0$ of the exterior algebra $\Omega(M,E)$.
That means that for all $\eta$ and $\zeta\in\Omega(M,E)$,
 \begin{equation*}{\cal L}_\rho(V)(\eta\wedge\zeta)=\bigl({\cal
 L}_\rho(V)\eta\bigr)\wedge\zeta+\eta\wedge\bigl({\cal
 L}_\rho(V)\zeta\bigr)\,.\end{equation*}
\par\smallskip
{\rm 4.}  For each $V$ and $W\in A^1(M,E)$, $\eta\in\Omega(M,E)$,
 \begin{equation*}{\cal L}_\rho\bigl(\{V,W\}\bigr)\eta=\bigl(
 {\cal L}_\rho(V)\circ {\cal L}_\rho(W)-{\cal L}_\rho(W)\circ
 {\cal L}_\rho(V)\bigr)\eta\,.\end{equation*}
\par\smallskip
{\rm 5.}  For each $V\in A^1(M,E)$, $f\in C^\infty(M,\RR)$ and
$\eta\in\Omega(M,E)$,
 \begin{equation*}{\cal L}_\rho(fV)\eta=f{\cal
 L}_\rho(V)\eta+d_\rho f\wedge i(V)\eta\,.\end{equation*}
\end{nprop}

\begin{nproof}
{\rm 1.} Let $W\in A^1(M,E)$. Then
 \begin{equation*}
 \begin{split}
 \bigl\langle{\cal L}_\rho(V)(d_\rho f),W\bigr\rangle
 &={\cal L}_\rho(V)\langle d_\rho f,W\rangle-\langle d_\rho
 f,\{V,W\}\rangle\\
 &={\cal L}(\rho\circ V)\circ{\cal L}(\rho\circ W)f
 -{\cal L}\bigl(\rho\circ\{V,W\}\bigr)f\\
 &={\cal L}(\rho\circ W)\circ{\cal L}(\rho\circ V)f\\
 &=\bigl\langle d_\rho\bigl({\cal
 L}_\rho(V)f\bigr),W\bigr\rangle\,,
 \end{split}
 \end{equation*}
so Property~1 is proven.
\par\medskip
{\rm 2.} Let $V$ and $W\in A^1(M,E)$, $\eta\in \Omega^p(M,E)$,
$V_1,\ldots,V_{p-1}\in A^1(M,E)$. We may write
 \begin{equation*}
 \begin{split}
 \bigl({\cal L}_\rho(V)\circ
 i(W)\eta\bigr)(V_1,
 &\ldots,V_{p-1})\\
 &={\cal L}_\rho(V)\bigl(\eta(W,V_1,\ldots,V_{p-1})\bigr)\\
 &\quad-\sum_{k=1}^{p-1}\eta\bigl(W,V_1,\ldots,V_{k-1},\{V,V_k\},
 V_{k+1},\ldots,V_{p-1}\bigr)\\
 &=\bigl({\cal L}_\rho(V)\eta\bigr)(W,V_1,\ldots,V_{p-1})
 +\eta\bigl(\{V,W\},V_1,\ldots, V_{p-1}\bigr)\\
 &=\Bigl(\bigl(i(W)\circ{\cal
 L}_\rho(V)+i(\{V,W\})\bigr)\eta\Bigr)
 (V_1,\ldots,V_{p-1})\,,
 \end{split}
 \end{equation*}
so Property~2 is proven.
\par\medskip
{\rm 3.} Let $V\in A^1(M,E)$, $\eta\in \Omega^p(M,E)$ and
$\zeta\in \Omega^q(M,E)$. For $p<0$, as well as for $q<0$, both
sides of the equality stated in Property~3 vanish, so that
equality is trivially satisfied. For $p=q=0$, that equality is
also satisfied, as shown by Equality $(*)$ in the proof of
Proposition~4.1.1. We still have to prove that equality for $p>0$
and (or) $q>0$. We will do that by induction on the value of
$p+q$. Let $r\geq 1$ be an integer such that the equality stated
in Property~3 holds for $p+q\leq r-1$. Such an integer exists, for
example $r=1$. We assume now that $p\geq 0$ and $q\geq 0$ are such
that $p+q=r$. Let $W\in A^1(M,E)$. By using Property~2, we may
write
 \begin{equation*}
 \begin{split}
 i(W)\circ{\cal L}_\rho(V)(\eta_\wedge\zeta)
 &={\cal L}_\rho(V)\circ
 i(W)(\eta\wedge\zeta)-i(\{V,W\})(\eta\wedge\zeta)\\
 &={\cal L}_\rho(V)\bigl(
 i(W)\eta\wedge\zeta+(-1)^p\eta\wedge i(W)\zeta\bigr)\\
 &\quad -i(\{V,W\})\eta\wedge\zeta-(-1)^p\eta\wedge
 i(\{V,W\})\zeta\,.
 \end{split}
 \end{equation*}
Since $i(W)\eta\in\Omega^{p-1}(M,E)$ and
$i(W)\zeta\in\Omega^{q-1}(M,E)$, the induction assumption allows
us to use Property~3 to transform the first terms of the right
hand side. We obtain
 \begin{equation*}
 \begin{split}
 i(W)\circ{\cal L}_\rho(V)(\eta_\wedge\zeta)
 &=\bigl({\cal L}_\rho(V)\circ
 i(W)\eta\bigr)\wedge\zeta+i(W)\eta\wedge{\cal L}_\rho(V)\zeta\\
 &\quad+(-1)^p\bigl({\cal L}_\rho(V)\eta\bigr)\wedge
 i(W)\zeta
 +(-1)^p\eta\wedge\bigl({\cal L}_\rho(V)\circ
 i(W)\zeta\bigr)\\
 &\quad-i(\{V,W\})\eta\wedge\zeta-(-1)^p\eta\wedge
 i(\{V,W\})\zeta\,.
 \end{split}
 \end{equation*}
By rearranging the terms, we obtain
 \begin{equation*}
 \begin{split}
 i(W)\circ{\cal L}_\rho(V)(\eta_\wedge\zeta)
 &=\bigl({\cal L}_\rho(V)\circ
 i(W)\eta-i(\{V,W\})\eta\bigr)\wedge\zeta\\
 &\quad+(-1)^p\eta\wedge\bigl({\cal L}_\rho(V)\circ i(W)\zeta
 -i(\{V,W\})\zeta\bigr)\\
 &\quad+i(W)\eta\wedge{\cal L}_\rho(V)\zeta
 +(-1)^p\bigl({\cal L}_\rho(V)\eta\bigr)\wedge i(W)\zeta\,.
 \end{split}
 \end{equation*}
By using again Property~2 we get
 \begin{equation*}
 \begin{split}
 i(W)\circ{\cal L}_\rho(V)(\eta_\wedge\zeta)
 &=\bigl(i(W)\circ{\cal L}_\rho(V)\eta\bigr)\wedge\zeta
   +(-1)^p\eta\wedge\bigl(i(W)\circ{\cal L}_\rho(V)\zeta\bigr)\\
 &\quad+i(W)\eta\wedge{\cal L}_\rho(V)\zeta
   +(-1)^p{\cal L}_\rho(V)\eta\wedge i(W)\zeta\\
 &=i(W)\bigl({\cal L}_\rho(V)\eta\wedge\zeta+\eta\wedge{\cal
 L}_\rho(V)\zeta\bigr)\,.
 \end{split}
 \end{equation*}
Since that last equality holds for all $W\in A^1(M,E)$, it follows
that Property~3 holds for $\eta\in\Omega^p(M,E)$ and
$\zeta\in\Omega^q(M,E)$, with $p\geq 0$, $q\geq 0$ and $p+q=r$. We
have thus proven by induction that Property~3 holds for all $p$
and $q\in\ZZ$, $\eta\in\Omega^p(M,E)$, $\zeta\in\Omega^q(M,E)$.
The same equality holds, by bilinearity, for all $\eta$ and
$\zeta\in\Omega(M,E)$.
\par\medskip
{\rm 4.} Let $V$ and $W\in A^1(M,E)$. Then $\{V,W\}\in A^1(M,E)$
and, by Property~3, ${\cal L}_\rho(V)$, ${\cal L}_\rho(W)$ and
${\cal L}_\rho(\{V,W\})$ are derivations of degree $0$ of the
graded algebra $\Omega(M,E)$. By 3.1.6~(ii), the graded bracket
 \begin{equation*}\bigl[{\cal L}_\rho(V),{\cal L}_\rho(W)\bigr]
 ={\cal L}_\rho(V)\circ {\cal L}_\rho(W)-{\cal L}_\rho(W)\circ
 {\cal L}_\rho(V)\end{equation*}
is also a derivation of degree $0$ of $\Omega(M,E)$. Property~4
means that the derivations ${\cal L}_\rho(\{V,W\})$ and
$\bigl[{\cal L}_\rho(V),{\cal L}_\rho(W)\bigr]$ are equal. In
order to prove that equality, it is enough to prove that it holds
true for $\eta\in \Omega^0(M,E)$ and for $\eta\in\Omega^1(M,E)$,
since the graded algebra $\Omega(M,E)$ is generated by its
homogeneous elements of degrees $0$ and $1$.
\par\smallskip
Let $f\in\Omega^0(M,E)=C^\infty(M,\RR)$. We have
 \begin{equation*}
 \begin{split}
 {\cal L}_\rho(\{V,W\})f
 &={\cal L}\bigl(\rho\circ\{V,W\}\bigr)f\\
 &={\cal L}\bigl([\rho\circ V,\rho\circ W]\bigr)f\\
 &=\bigl[{\cal L}(\rho\circ V),{\cal L}(\rho\circ W)\bigr]f\\
 &=\bigl[{\cal L}_\rho(V),{\cal L}_\rho(W)\bigr]f\,,
 \end{split}
 \end{equation*}
therefore Property~4 holds for $\eta=f\in\Omega^0(M,E)$.
\par\smallskip
Now let $\eta\in\Omega^1(M,E)$ and $Z\in A^1(M,E)$. By using
Property~2, then Property~4 for elements $\eta\in\Omega^0(M,E)$,
we may write
\goodbreak
 \begin{equation*}
 \begin{split}
 i(Z)\circ{\cal L}_\rho(\{V,W\})\eta
 &={\cal
 L}_\rho(\{V,W\})\bigl(i(Z)\eta\bigr)-i\bigl(\bigl\{\{V,W\},Z\bigr\}
 \bigr)\eta\\
 &=\bigl({\cal L}_\rho(V)\circ {\cal L}_\rho(W)
   -{\cal L}_\rho(W)\circ {\cal
   L}_\rho(V)\bigr)\bigl(i(Z)\eta\bigr)\\
  &\quad-i\bigl(\bigl\{\{V,W\},Z\bigr\}\bigr)\eta\,.
 \end{split}
  \end{equation*}
By using Property~2 and the Jacobi identity, we obtain
 \begin{equation*}
 \begin{split}
 i(Z)\circ{\cal L}_\rho(\{V,W\})\eta
 &={\cal L}_\rho(V)\bigl(i(\{W,Z\})\eta+i(Z)\circ{\cal
 L}_\rho(W)\eta\bigr)\\
 &\quad-{\cal L}_\rho(W)\bigl(i(\{V,Z\})\eta+i(Z)\circ{\cal
 L}_\rho(V)\eta\bigr)\\
 &\quad-i\bigl(\bigl\{\{V,W\},Z\bigr\}\bigr)\eta\\
 &=i(\{W,Z\}){\cal L}_\rho(V)\eta + i(\{V,Z\}){\cal
 L}_\rho(W)\eta\\
 &\quad-i(\{V,Z\}){\cal L}_\rho(W)\eta
 - i(\{W,Z\}){\cal L}_\rho(V)\eta\\
 &\quad+i(Z)\circ\bigl({\cal L}_\rho(V)\circ{\cal L}_\rho(W)
 -{\cal L}_\rho(W)\circ{\cal L}_\rho(V)\bigr)\eta\\
 &\quad+i\Bigl(\bigl\{V,\{W,Z\}\bigr\}-\bigl\{W,\{V,Z\}\bigr\}
 -\bigl\{\{V,W\},Z\bigr\}\Bigr)\eta\\
 &=i(Z)\circ\bigl({\cal L}_\rho(V)\circ{\cal L}_\rho(W)
 -{\cal L}_\rho(W)\circ{\cal L}_\rho(V)\bigr)\eta\,.
 \end{split}
 \end{equation*}
Since that last equality holds for all $Z\in A^1(M,E)$, Property~4
holds for all $\eta\in\Omega^1(M,E)$, and therefore for all
$\eta\in\Omega(M,E)$.
\par\medskip
{\rm 5.} Let $V\in A^1(M,E)$ and $f\in C^\infty(M,\RR)$. We have
seen (Property~4) that ${\cal L}_\rho(fV)$ is a derivation of
degree $0$ of $\Omega(M,E)$. We easily verify that
 \begin{equation*}\eta\mapsto f{\cal L}_\rho(V)\eta+d_\rho f\wedge
i(V)\eta\end{equation*} is too a derivation of degree $0$ of
$\Omega(M,E)$. Property~5 means that these two derivations are
equal. As above, it is enough to prove that Property~5 holds for
$\eta\in\Omega^0(M,E)$ and for $\eta\in\Omega^1(M,E)$.
\par\smallskip
Let $g\in\Omega^0(M,E)=C^\infty(M,\RR)$. We may write
 \begin{equation*}{\cal L}_\rho(fV)g=i(fV)d_\rho g=f{\cal
L}_\rho(V)g\,,\end{equation*} which shows that Property~5 holds
for $\eta=g\in\Omega^0(M,E)$.
\par\smallskip
Let $\eta\in\Omega^1(M,E)$, and $W\in A^1(M,E)$. We have
 \begin{equation*}
 \begin{split}
 \bigl\langle {\cal L}_\rho(fV)\eta,W\bigr\rangle
 &={\cal L}_\rho(fV)\bigl(\langle\eta,W\rangle\bigr)
  -\bigl\langle\eta,\{fV,W\}\bigr\rangle\\
 &=f{\cal L}_\rho(V)\bigl(\langle\eta,W\rangle\bigr)
  -f\bigl\langle\eta,\{V,W\}\bigr\rangle
  +\bigl\langle\eta,\bigl(i(W)d_\rho
  f\bigr)V\bigr\rangle\\
 &=\bigl\langle f{\cal L}_\rho(V)\eta,W\bigr\rangle
  +\bigl(i(W)d_\rho f\bigr)i(V)\eta\\
 &=\bigl\langle f{\cal L}_\rho(V)\eta +d_\rho f\wedge i(V)\eta,W\bigr\rangle
 \,,
 \end{split}
 \end{equation*}
since, $\eta$ being in $\Omega^1(M,E)$, $i(W)\circ i(V)\eta=0$.
The last equality being satisfied for all $W\in A^1(M,E)$, the
result follows.
\end{nproof}
\par\bigskip
The next Proposition shows that for each $V\in A^1(M,E)$, the Lie
derivative ${\cal L}_\rho(V)$, already defined as a derivation of
degree $0$ of the graded algebra $\Omega(M,E)$, can also be
extended into a derivation of degree $0$ of the graded algebra
$A(M,E)$, with very nice properties. As we will soon see, the
Schouten-Nijenhuis bracket will appear as a very natural further
extension of the Lie derivative.

\begin{nprop}\sl Let $(E,\tau,M,\rho)$ be a Lie algebroid on a
smooth manifold $M$. For each smooth section $V\in A^1(M,E)$ of
the vector bundle $(E,\tau,M)$, there exists a unique graded
endomorphism of degree $0$ of the graded algebra of multivectors
$A(M,E)$, called the {\it Lie derivative\/} with respect to $V$
and denoted by ${\cal L}_\rho(V)$, which satisfies the following
properties:

\begin{description}
\item{\rm(i)} For a smooth function
$f\in A^0(M,E)=C^\infty(M,\RR)$,
 \begin{equation*}{\cal L}_\rho(V)f=i(\rho\circ V)df={\cal L}(\rho\circ
V)f\,,\end{equation*} where ${\cal L}(\rho\circ V)$ denotes the
usual Lie derivative with respect to the vector field $\rho\circ
V$;

\item{\rm (ii)} For an integer $p\geq 1$ and a multivector
$P\in A^p(M,E)$, ${\cal L}_\rho(V)P$ is the unique element in
$A^p(M,E)$ such that, for all $\eta\in \Omega^p(M,E)$,
 \begin{equation*}\bigl\langle\eta,{\cal L}_\rho(V)P\bigr\rangle
 ={\cal L}_\rho(V)\bigl(\langle\eta,P\rangle\bigr)-\bigl\langle
 {\cal L}_\rho(V)\eta,P\bigr\rangle\,.\end{equation*}
\end{description}
\end{nprop}

\begin{nproof} Let us first observe that
$A^0(M,E)=\Omega^0(M,E)=C^\infty(M,\RR)$, and that for $f\in
A^0(M,E)$, the definition of ${\cal L}_\rho(V)f$ given above is
the same as that given in Proposition 4.1.1.
\par\smallskip
Now let $p\geq 1$ and $P\in A^p(M,E)$. The map
 \begin{equation*}\eta\mapsto K(\eta)=
 {\cal L}_\rho(V)\bigl(\langle\eta,P\rangle\bigr)-\bigl\langle
 {\cal L}_\rho(V)\eta,P\bigr\rangle\end{equation*}
is clearly an $\RR$-linear map defined on $\Omega^p(M,E)$, with
values in $C^\infty(M,\RR)$. Let $f\in C^\infty(M,\RR)$. We have
 \begin{equation*}
 \begin{split}
 K(f\eta)
 &={\cal L}_\rho(V)\bigl(\langle f\eta,P\rangle\bigr)-\bigl\langle
 {\cal L}_\rho(V)(f\eta),P\bigr\rangle\\
 &=f\Bigl({\cal L}_\rho(V)\bigl(\langle\eta,P\rangle\bigr)-\bigl\langle
 {\cal L}_\rho(V)\eta,P\bigr\rangle\Bigr)\\
 &\quad+\bigl({\cal L}_\rho(V)f\bigr)\langle\eta,P\rangle
        -\bigl({\cal L}_\rho(V)f\bigr)\langle\eta,P\rangle\\
 &=fK(\eta)\,.
 \end{split}
 \end{equation*}
This proves that the map $K$ is $C^\infty(M,\RR)$-linear. Since
the pairing allows us to consider the vector bundle
$(\bigwedge^pE,\tau,M)$ as the dual of $(\bigwedge^pE^*,\pi,M)$,
we see that there exists a unique element ${\cal L}_\rho(V)P\in
A^p(M,E)$ such that, for all $\eta\in \Omega^p(M,E)$,
 \begin{equation*}K(\eta)=
 {\cal L}_\rho(V)\bigl(\langle\eta,P\rangle\bigr)-\bigl\langle
 {\cal L}_\rho(V)\eta,P\bigr\rangle=\bigl\langle\eta,{\cal
 L}_\rho(V)P\bigr\rangle\,,\end{equation*}
and that ends the proof.
\end{nproof}

\begin{nprop}\sl Under the assumptions of Proposition 4.1.5, the
Lie derivative has the following properties.
\par\smallskip
{\rm 1.} For each $V$ and $W\in A^1(M,E)$,
 \begin{equation*}{\cal L}_\rho(V)(W)=\{V,W\}\,.\end{equation*}
\par\smallskip
{\rm 2.} For $V,V_1,\ldots,V_p\in A^1(M,E)$,
 \begin{equation*}{\cal L}_\rho(V)(V_1\wedge\cdots\wedge V_p)
 =\sum_{i=1}^pV_1\wedge\cdots\wedge V_{i-1}\wedge\{V,V_i\}\wedge V_{i+1}
 \wedge\cdots\wedge V_p\,.\end{equation*}
\par\smallskip
{\rm 3.} For each $V\in A^1(M,E)$, ${\cal L}_\rho(V)$ is a
derivation of degree $0$ of the exterior algebra $A(M,E)$. That
means that for all $P$ and $Q\in A(M,E)$,
 \begin{equation*}{\cal L}_\rho(V)(P\wedge Q)=
   \bigl({\cal L}_\rho(V)P\bigr)\wedge Q+
   P\wedge{\cal L}_\rho(V)Q\,.\end{equation*}
\par\smallskip
{\rm 4.} For each $V\in A^1(M,E)$, $P\in A(M,E)$ and
$\eta\in\Omega(M,E)$,
 \begin{equation*}{\cal L}_\rho(V)\bigl(i(P)\eta\bigr)=i\bigl({\cal
 L}_\rho(V)P\bigr)\eta
 +i(P)\bigl({\cal L}_\rho(V)\eta\bigr)\,.\end{equation*}
\par\smallskip
{\rm 5.} Similarly, for each $V\in A^1(M,E)$, $P\in A(M,E)$ and
$\eta\in\Omega(M,E)$,
 \begin{equation*}{\cal L}_\rho(V)\bigl(\langle\eta,P\rangle\bigr)
 =\bigl\langle{\cal L}_\rho(V)\eta,P\bigr\rangle
 +\bigl\langle\eta,{\cal L}_\rho(V)P\bigr\rangle\,.\end{equation*}
\par\smallskip
{\rm 6.}  For each $V$ and $W\in A^1(M,E)$, $P\in A(M,E)$,
 \begin{equation*}{\cal L}_\rho\bigl(\{V,W\}\bigr)P=\bigl(
 {\cal L}_\rho(V)\circ {\cal L}_\rho(W)-{\cal L}_\rho(W)\circ
 {\cal L}_\rho(V)\bigr)P\,.\end{equation*}
\par\smallskip
{\rm 7.}  For each $V\in A^1(M,E)$, $f\in C^\infty(M,\RR)$, $P\in
A(M,E)$ and $\eta\in \Omega(M,E)$,
 \begin{equation*}\bigl\langle\eta,{\cal L}_\rho(fV)P\bigr\rangle
 =f\bigl\langle{\cal L}_\rho(V)P,\eta\bigr\rangle
 +\bigl\langle d_\rho f\wedge i(V)\eta,P\bigr\rangle\,.\end{equation*}
\end{nprop}

\begin{nproof}
{\rm 1.} Let $V$ and $W\in A^1(M,E)$, $\eta\in\Omega(M,E)$. We may
write
 \begin{equation*}
 \begin{split}
 \bigl\langle \eta,{\cal L}_\rho(V)W\bigr\rangle
 &={\cal L}_\rho(V)\bigl(\langle\eta,W\rangle\bigr)-\bigl\langle{\cal
 L}_\rho(V)\eta,W\bigr\rangle\\
 &=\bigl\langle\eta,\{V,W\}\bigr\rangle\,.
 \end{split}
 \end{equation*}
We have proven Property~1.
\par\smallskip
{\rm 2.} The proof is similar to that of Property~1.
\par\smallskip
{\rm 3.} When $P=V_1\wedge\cdots\wedge V_p$ and $Q=W_1\wedge\cdots
\wedge W_q$ are decomposable homogeneous elements in $A(M,E)$,
Property~3 is an easy consequence of 2. The validity of Property~3
for all $P$ and $Q\in A(M,E)$ follows by linearity.
\par\smallskip
{\rm 4.}  When $P=V_1\wedge\cdots\wedge V_p$ is a decomposable
homogeneous element in $A(M,E)$, Property~4 is an easy consequence
of Property~2. The validity of Property~4 for all $P$ and $Q\in
A(M,E)$ follows by linearity.
\par\smallskip
{\rm 5.} This is an immediate consequence of Property~4.
\par\smallskip
{\rm 6.} This is an immediate consequence of Property~4 of this
Proposition and of Property~4 of Proposition 4.1.4.
\par\smallskip
{\rm 7.} This is an immediate consequence of Property~4 of this
Proposition and of Property~5 of Proposition 4.1.4.
\end{nproof}

\subsection{The $\Omega(M,E)$-valued exterior derivative}

We have introduced above (Definition 4.1.2) the
$\Omega(M,E)$-valued exterior derivative of a function $f\in
\Omega^0(M,E)=C^\infty(M,\RR)$. The next proposition shows that
the $\Omega(M,E)$-valued exterior derivative extends as a graded
endomorphism of degree $1$ of the graded algebra $\Omega(M,E)$. We
will see later (Proposition 4.2.3) that the $\Omega(M,E)$-valued
exterior derivative is in fact a derivation of degree $1$ of
$\Omega(M,E)$.

\begin{nprop}\sl Let $(E,\tau,M,\rho)$ be a Lie algebroid over a
smooth manifold $M$. There exists a unique graded endomorphism of
degree $1$ of the exterior algebra of forms $\Omega(M,E)$, called
the {\it $\Omega(M,E)$-valued exterior derivative\/} (or, in
brief, the {\it exterior derivative\/}) and denoted by $d_\rho$,
which satisfies the following properties:
\par\smallskip
\item{\rm(i)} For $f\in\Omega^0(M,E)=C^\infty(M,\RR)$, $d_\rho f$ is the
unique element in $\Omega^1(M,E)$, already defined (Definition
4.1.2), such that, for each $V\in A^1(M,E)$,
 \begin{equation*}\langle d_\rho f,V\rangle={\cal L}_\rho(V)f=\langle
df,\rho\circ
 V\rangle=\langle{}^t\!\rho\circ df,V\rangle\,,\end{equation*}
where $d$ stands for the usual exterior derivative of smooth
functions on $M$, and ${}^t\!\rho:E^*\to T^*M$ is the transpose of
the anchor $\rho$.
\par\smallskip
\item{\rm(ii)} For $p\geq 1$ and $\eta\in \Omega^p(M,E)$, $d_\rho\eta$ is
the unique element in $\Omega^{p+1}(M,E)$ such that, for all
$V_0,\ldots, V_p\in A^1(M,E)$,
 \begin{equation*}
 \begin{split}
 d_\rho\eta(V_0,\ldots,V_p)
 &=\sum_{i=0}^p(-1)^i{\cal L}_\rho(V_i)\bigl(\eta(V_0,\ldots,\widehat{V}_i,
 \ldots,V_p)\bigr)\\
 &\quad+\sum_{0\leq i<j\leq p}(-1)^{i+j}\eta\bigl(\{V_i,V_j\},V_0,\ldots,
 \widehat{V}_i,\ldots,\widehat{V}_j,\ldots,V_p)\bigr)\,,
 \end{split}
 \end{equation*}
where the symbol~~$\widehat{\ }$~~over the terms $V_i$ and $V_j$
means that these terms are omitted.

\end{nprop}

\begin{nproof} For $f\in\Omega^0(M,E)$, $d_\rho f$ is clearly an element
in $\Omega^1(M,E)$.
\par
Let $p\geq 1$ and $\eta\in\Omega^p(M,E)$. As defined in {\rm(ii)},
$d_\rho\eta$ is a map, defined on $\bigl(A^1(M,E)\bigr)^p$, with
values in $C^\infty(M,\RR)$. The reader will immediately see that
this map is skew-symmetric and $\RR$-linear in each of its
arguments. In order to prove that $d_\rho\eta$ is an element in
$\Omega^{p+1}(M,E)$, it remains only to verify that as a map,
$d_\rho\eta$ is $C^\infty(M,\RR)$-linear in each of its arguments,
or simply in its first argument, since the skew-symmetry will
imply the same property for all other arguments. Let $f\in
C^\infty(M,\RR)$. We have
 \begin{equation*}
 \begin{split}
 d_\rho\eta(fV_0,&V_1,\ldots,V_p)\\
 &={\cal L}_\rho(fV_0)\bigl(\eta(V_1,\ldots,V_p)\bigr)\\
 &\quad+\sum_{i=1}^p(-1)^i{\cal
 L}_\rho(V_i)\bigl(f\eta(V_0,\ldots,\widehat{V}_i,\ldots,V_p)\bigr)\\
 &\quad+\sum_{1\leq j\leq
 p}(-1)^j\eta\bigl(\{fV_0,V_j\},V_1,\ldots,\widehat{V}_j,\ldots,V_p\bigr)\\
 &\quad+\sum_{1\leq i<j\leq p}(-1)^{i+j}\eta\bigl(\{V_i,V_j\},fV_0,V_1,
 \ldots,\widehat{V}_i,\ldots,\widehat{V}_j,\ldots,V_p\bigr)\,.
 \end{split}
 \end{equation*}
By a rearrangement of the terms in the right hand side, and by
using the formulae
 \begin{equation*}{\cal L}_\rho(V_i)\bigl(f\eta(\ldots)\bigr)=\bigl({\cal
 L}_\rho(V_i)f\bigr)\eta(\ldots) + f{\cal
 L}_\rho(V_i)\bigl(\eta(\ldots)\bigr)\,,\end{equation*}
and
 \begin{equation*}\{fV_0,V_j\}=f\{V_0,V_j\}-\bigl({\cal
L}_\rho(V_j)f\bigr)V_0\,,\end{equation*} we obtain
 \begin{equation*}d_\rho\eta(fV_0,V_1,\ldots,V_p)=fd_\rho\eta(V_0,V_1,\ldots,V_p
)\,.\end{equation*} We have shown that
$d_\rho\eta\in\Omega^{p+1}(M,E)$. The $\Omega(M,E)$-valued
exterior derivative so defined on $\Omega^p(M,E)$ for all
$p\in\ZZ$ extends, in a unique way, into a graded endomorphism of
degree $1$ of $\Omega(M,E)$. \end{nproof}

\begin{nrmk}\rm Let $p\geq 1$, $\eta\in\Omega^p(M,E)$ and
$V_0,\ldots,V_p\in A^1(M,E)$. The formula for $d_\rho\eta$ given
in Proposition 4.2.1 can be cast into another form, often useful:
 \begin{equation*}
 \begin{split}
 d_\rho\eta(V_0,\ldots,V_p)
 &=\sum_{i=0}^p(-1)^i\bigl({\cal
 L}_\rho(V_i)\eta\bigr)(V_0,\ldots,\widehat{V_i},
 \ldots,V_p)\\
 &\quad-\sum_{0\leq i<j\leq p}(-1)^{i+j}\eta\bigl(\{V_i,V_j\},V_0,\ldots,
 \widehat{V}_i,\ldots,\widehat{V}_j,\ldots,V_p)\bigr)\,,
 \end{split}
 \end{equation*}
where the symbol~~$\widehat{\ }$~~over the terms $V_i$ and $V_j$
means that these terms are omitted.
\par\smallskip
For example, for $p=1$,
 \begin{equation*}
 \begin{split}
 d_\rho\eta(V_0,V_1)
 &={\cal L}_\rho(V_0)\bigl(\eta(V_1)\bigr)
   -{\cal L}_\rho(V_1)\bigl(\eta(V_0)\bigr)
   -\eta\bigl(\{V_0,V_1\}\bigr)\\
 &=\bigl\langle{\cal L}_\rho(V_0)\eta,V_1\bigr\rangle
   -\bigl\langle{\cal L}_\rho(V_1)\eta,V_0\bigr\rangle
   +\eta\bigl(\{V_0,V_1\}\bigr)\,.
 \end{split}
   \end{equation*}
\end{nrmk}

\begin{nprop}\sl Under the assumptions of Proposition 4.2.1,
the $\Omega(M,E)$-valued exterior derivative has the following
properties:
\par\smallskip
{\rm 1.} Let $V\in A^1(M,E)$. The Lie derivative ${\cal
L}_\rho(V)$, the exterior derivative $d_\rho$ and the interior
product $i(V)$ are related by the formula
 \begin{equation*}{\cal L}_\rho(V)=i(V)\circ d_\rho+d_\rho\circ
i(V)\,.\end{equation*}
\par\smallskip
{\rm 2.} The exterior derivative $d_\rho$ is a derivation of
degree $1$ of the exterior algebra $\Omega(M,E)$. That means that
for each $p\in\ZZ$, $\eta\in\Omega^p(M,E)$ and
$\zeta\in\Omega(M,E)$,
 \begin{equation*}d_\rho(\eta\wedge\zeta)=d_\rho\eta\wedge\zeta+(-1)^p\eta\wedge
 d_\rho\zeta\,.\end{equation*}
\par\smallskip
{\rm 3.} Let $V\in A^1(M,E)$. Then
 \begin{equation*}{\cal L}_\rho(V)\circ d_\rho=d_\rho\circ{\cal
L}_\rho(V)\,.\end{equation*}
\par\smallskip
{\rm 4.} The square of $d_\rho$ vanishes identically,
 \begin{equation*}d_\rho\circ d_\rho=0\,.\end{equation*}

\end{nprop}

\begin{nproof}
{\rm 1.} Let $V_0=V,\ V_1,\ \ldots,\ V_p\in A^1(M,E)$,
$\eta\in\Omega^p(M,E)$. Then
 \begin{equation*}
 \begin{split}
 \bigl(i(V)\circ d_\rho \eta\bigr)(V_1,
 &\ldots,V_p)\\
 &=d_\rho\eta(V,V_1,\ldots,V_p)\\
 &=\sum_{i=0}^p(-1)^i{\cal
 L}_\rho(V_i)\bigl(\eta(V_0,\ldots,\widehat{V}_i,\ldots,V_p)\bigr)\\
 &\quad+\sum_{0\leq i<j\leq p}(-1)^{i+j}\eta\bigl(\{V_i,V_j\},V_0,\ldots,
 \widehat{V}_i,\ldots,\widehat{V}_j,\ldots,V_p\bigr)\,,
 \end{split}
 \end{equation*}
and
 \begin{equation*}
 \begin{split}
 \bigl(d_\rho\circ i(V) \eta\bigr)(V_1,
 &\ldots,V_p)\\
 &=\sum_{i=1}^p(-1)^{i-1}{\cal
 L}_\rho(V_i)\bigl(\eta(V_0,\ldots,\widehat{V_i},\ldots,V_p)\bigr)\\
 &\quad+\sum_{1\leq i<j\leq
 p}(-1)^{i+j}\eta\bigl(V_0,\{V_i,V_j\},V_1,\ldots,
 \widehat{V}_i,\ldots,\widehat{V}_j,\ldots,V_p\bigr)\,.
 \end{split}
 \end{equation*}
Let us add these two equalities. Several terms cancel, and we
obtain
 \begin{equation*}
 \begin{split}
 \Bigl(\bigl(i(V)\circ d_\rho
 &+d_\rho\circ i(V)\bigr)\eta\Bigr)(V_1,\ldots,V_p)\\
 &={\cal L}_\rho(V_0)\bigl(\eta(V_1,\ldots,V_p)\bigr)
 +\sum_{j=1}^p(-1)^j\eta\bigl(\{V_0,V_j\},V_1,
 \ldots,\widehat{V}_j,\ldots,
 V_p\bigr)\,.
 \end{split}
 \end{equation*}
When we shift, in the last term of the right hand side, the
argument $\{V_0,V_j\}$ to the $j$-th position, we obtain
 \begin{equation*}
 \begin{split}
 \Bigl(\bigl(i(V)&\circ d_\rho+d_\rho\circ i(V)\bigr)
 \eta\Bigr)(V_1,\ldots,V_p)\\
 &={\cal L}_\rho(V_0)\bigl(\eta(V_1,\ldots,V_p)\bigr)
 +\sum_{j=1}^p\eta\bigl(V_1,\ldots,V_{j-1},\{V_0,V_j\},V_{j+1},\ldots,
 V_p\bigr)\\
 &=\bigl({\cal L}_\rho(V_0)\eta\bigr)(V_1,\ldots,V_p)\,.
 \end{split}
 \end{equation*}
\par\smallskip
{\rm 2.} For $\eta=f$ and
$\zeta=g\in\Omega^0(M,E)=C^\infty(M,\RR)$, Property 2 holds since
we have, for alll $V\in A^1(M,E)$,
 \begin{equation*}
 \begin{split}
 \bigl\langle d_\rho(fg),V\bigr\rangle
 =\bigl\langle d(fg),\rho\circ V\bigr\rangle
 &=\langle f\,dg+g\,df,\rho\circ V\rangle\\
 &=\langle f\,d_\rho g+g\,d_\rho f,V\rangle\,.
 \end{split}
 \end{equation*}
Now let $p\geq 0$ and $q\geq 0$ be two integers, $\eta\in
\Omega^p(M,E)$, $\zeta\in\Omega^q(M,E)$. We will prove that
Property 2 holds by induction on $p+q$. Just above, we have seen
that it holds for $p+q=0$. Let us assume that $r$ is an integer
such that Property 2 holds for $p+q\leq r$, and that now
$p+q=r+1$. Let $V\in A^1(M,E)$. We may write
 \begin{equation*}
 \begin{split}
 i(V)d_\rho(\eta\wedge\zeta)
 &={\cal L}_\rho(V)(\eta\wedge\zeta)-d_\rho\circ
 i(V)(\eta\wedge\zeta)\\
 &={\cal L}_\rho(V)\eta\wedge\zeta+\eta\wedge{\cal
 L}_\rho(V)\zeta\\
 &\quad-d_\rho\bigl(i(V)\eta\wedge\zeta+(-1)^p\eta\wedge
 i(V)\zeta\bigr)\,.
 \end{split}
 \end{equation*}
We may now use the induction assumption, since in the last terms
of the right hand side $i(V)\eta\in\Omega^{p-1}(M,E)$ and
$i(V)\zeta\in\Omega^{q-1}(M,E)$. After some rearrangements of the
terms we obtain
 \begin{equation*}i(V)d_\rho(\eta\wedge\zeta)
 =i(V)(d_\rho\eta\wedge\zeta+\eta\wedge
 d_\rho\zeta)\,.
 \end{equation*}
Since that result holds for all $V\in A^1(M,E)$, Property 2 holds
for $p+q=r+1$, and therefore for all $p$ and $q\in\ZZ$.
\par\smallskip
{\rm 3.} Let $V\in A^1(M,E)$. We know (Proposition 4.1.4) that
${\cal L}_\rho(V)$ is a derivation of degree $0$ of the exterior
algebra $\Omega(M,E)$, and we have just seen (Property 2) that
$d_\rho$ is a derivation of degree $1$ of that algebra. Therefore,
by 3.1.6, their graded bracket
 \begin{equation*}\bigl[{\cal L}_\rho(V),d_\rho\bigr]={\cal L}_\rho(V)\circ
d_\rho
 -d_\rho\circ{\cal L}_\rho(V)\end{equation*}
is a derivation of degree $1$ of $\Omega(M,E)$. In order to prove
that that derivation is equal to $0$, it is enough to prove that
it vanishes on $\Omega^0(M,E)$ and on $\Omega^1(M,E)$. We have
already proven that it vanishes on $\Omega^0(M,E)$ (Property~1 of
4.1.4). Let $\eta\in\Omega^1(M,E)$ and $W\in A^1(M,E)$. By using
Property~1 of this Proposition and Property~2 of~4.1.4, we may
write
 \begin{equation*}
 \begin{split}
 i(W)\circ\bigl({\cal L}_\rho(V)\circ d_\rho
 &-d_\rho\circ {\cal
 L}_\rho(V)\bigr)\eta\\
 &={\cal L}_\rho(V)\circ i(W)\circ d_\rho\eta-i\bigl(\{V,W\}\bigr)\circ
 d_\rho\eta\\
 &\quad-{\cal L}_\rho(W)\circ {\cal L}_\rho(V)\eta
 +d_\rho\circ i(W)\circ {\cal L}_\rho(V)\eta\,.
 \end{split}
 \end{equation*}
By rearrangement of the terms, we obtain
 \begin{equation*}
 \begin{split}
 i(W)\circ\bigl({\cal L}_\rho(V)\circ d_\rho
 &-d_\rho\circ {\cal
 L}_\rho(V)\bigr)\eta\\
 &=\bigl({\cal L}_\rho(V)\circ{\cal L}_\rho(W)-{\cal L}_\rho(W)\circ
 {\cal L}_\rho(V)-{\cal L}_\rho(\{V,W\})\bigr)\eta\\
 &\quad+d_\rho\circ i(\{V,W\})\eta-d_\rho\circ i(\{V,W\})\eta\\
 &\quad-\bigl({\cal L}_\rho(V)\circ d_\rho
 -d_\rho\circ{\cal L}_\rho(V)\bigr)\bigl(i(W)\eta\bigr)\\
 &=0\,,
 \end{split}
 \end{equation*}
since $i(W)\eta\in\Omega^0(M,E)$, which implies that the last term
vanishes.
\par\smallskip
{\rm 4.} Property 2 shows that $d_\rho$ is a derivation of degree
$1$ of $\Omega(M,E)$. We know (3.1.6) that
$[d_\rho,d_\rho]=2d_\rho\circ d_\rho$ is a derivation of degree
$2$ of $\Omega(M,\RR)$. In order to prove that $d_\rho\circ
d_\rho=0$, it is enough to prove that it vanishes on
$\Omega^0(M,E)$ and on $\Omega^1(M,E)$.
\par
Let $f\in\Omega^0(M,E)=C^\infty(M,\RR)$, $V$ and $W\in A^1(M,E)$.
Then
 \begin{equation*}
 \begin{split}
 (d_\rho\circ d_\rho f)(V,W)
 &={\cal L}_\rho(V)\bigl(d_\rho f(W)\bigr)
 -{\cal L}_\rho(W)\bigl(d_\rho f(V)\bigr)
 -d_\rho f\bigl(\{V,W\}\bigr)\\
 &=\Bigl({\cal L}_\rho(V)\circ{\cal L}_\rho(W)
   -{\cal L}_\rho(W)\circ{\cal L}_\rho(V)
   -{\cal L}_\rho\bigl(\{V,W\}\bigr)\Bigr)f\\
 &=0\,,
 \end{split}
 \end{equation*}
where we have used Property~4 of Proposition 4.1.4. We have shown
that $d_\rho\circ d_\rho$ vanishes on $\Omega^0(M,E)$.
\par
Now let $\eta\in \Omega^1(M,E)$, $V_0$, $V_1$ and $V_2\in
A^1(M,E)$. Using Property 1, we may write
 \begin{equation*}
 \begin{split}
 (d_\rho\circ d_\rho \eta)(V_0,V_1,V_2)
 &=\Bigl(\bigl(i(V_0)\circ
 d_\rho\bigr)(d_\rho\eta)\Bigr)(V_1,V_2)\\
 &=\Bigl(\bigl({\cal L}_\rho(V_0)\circ d_\rho-d_\rho\circ
 i(V_0)\circ d_\rho\bigr)\eta\Bigr)(V_1,V_2)\,.
 \end{split}
 \end{equation*}
The last term in the right hand side may be transformed, by using
again Property 1:
 \begin{equation*}
 \begin{split}
 d_\rho\circ i(V_0)\circ d_\rho (\eta)
 &=d_\rho\circ {\cal L}_\rho(V_0)
 \eta-d_\rho\circ d_\rho\bigl(i(V_0)\eta\bigr)\\
 &=d_\rho\circ {\cal L}_\rho(V_0)\eta\,,
 \end{split}
 \end{equation*}
since, as $i(V_0)\eta\in \Omega^0(M,E)$, we have $d_\rho\circ
d_\rho\bigl(i(V_0)\eta\bigr)=0$. So we obtain
 \begin{equation*}(d_\rho\circ d_\rho \eta)(V_0,V_1,V_2)=
 \Bigl(\bigl({\cal L}_\rho(V_0)\circ d_\rho-d_\rho\circ
 {\cal L}_\rho(V_0)\bigr)\eta\Bigr)(V_1,V_2)\,.\end{equation*}
But Property 3 shows that
 \begin{equation*}\bigl({\cal L}_\rho(V_0)\circ d_\rho-d_\rho\circ
 {\cal L}_\rho(V_0)\bigr)\eta=0\,,\end{equation*}
so we have
 \begin{equation*}(d_\rho\circ d_\rho \eta)(V_0,V_1,V_2)=0\,,\end{equation*}
and our proof is complete.
\end{nproof}

\subsection{Defining a Lie algebroid by properties of its dual}
Let $(E,\tau,M)$ be a vector bundle and $(E^*,\pi,M)$ its dual
bundle. We have seen in 4.2 that when $(E,\tau,M)$ has a Lie
algebroid structure whose anchor is denoted by $\rho$, this
structure determines, on the graded algebra $\Omega(M,E)$ of
sections of the exterior powers of the dual bundle
$(E^*,\pi,M)$, a graded derivation $d_\rho$, of degree $1$,
which satisfies $d_\rho^2=d_\rho\circ d_\rho=0$. Now we are going
to prove a converse of this property: when a graded derivation of
degree $1$, whose square vanishes, is given on $\Omega(M,E)$, it
determines a Lie algebroid structure on $(E,\tau,M)$. This property
will be used later to prove that the cotangent bundle of a Poisson
manifold has a natural Lie algebroid structure.
\par\smallskip
We will need the following lemmas.

\begin{nlemma}\sl
Let $(E,\tau,M)$ be a vector bundle and $(E^*,\pi,M)$ its dual
bundle. Let $\delta$ be a graded derivation of degree $1$ of the
exterior algebra $\Omega(M,E)$ (notations defined in 3.2.4). For
each pair $(X,Y)$ of smooth sections of $(E,\tau,M)$ there exists a
unique smooth section $[X,Y]_\delta$ of $(E,\tau,M)$, called the
{\it $\delta$-bracket\/} of $X$ and $Y$, such that
 \begin{equation*}i\bigl([X,Y]_\delta\bigr)
 =\bigl[[i(X),\delta],i(Y)\bigr]\,.
 \end{equation*}
\end{nlemma}

\begin{nproof}
The map defined by the right hand side of the above equality,
  \begin{equation*}D:\eta\mapsto
D(\eta)=\bigl[[i(X),\delta],i(Y)\bigr]\end{equation*} is a
derivation of degree $-1$ of $\Omega(M,E)$, since it is obtained
by repeated application of the graded bracket to derivations
(property 3.1.6~(ii)). Therefore, it vanishes on
$\Omega^0(M,E)=C^\infty(M,\RR)$. As a consequence, that map is
$C^\infty(M,\RR)$-linear; we have indeed, for each $f\in
C^\infty(M,\RR)$ and $\eta\in\Omega(M,R)$,
 \begin{equation*}D(f\eta)=D(f)\wedge\eta+fD(\eta)=fD(\eta)\,.\end{equation*}
Therefore, there exists a unique smooth section $[X,Y]_\delta$ of
$(E,\tau,M)$ such that, for each $\eta\in \Omega^1(M,E)$,
 \begin{equation*}\bigl\langle
\eta,[X,Y]_\delta\bigr\rangle=D(\eta)\,.\end{equation*} Now the
maps
 \begin{equation*}i\bigl([X,Y]_\delta\bigr)\quad\hbox{and}\quad
   \bigl[[i(X),\delta],i(Y)\bigr]\end{equation*}
are both derivations of degree $-1$ of $\Omega(M,E)$, which
coincide on $\Omega^0(M,E)$ and $\Omega^1(M,E)$. Since derivations
are local, and since $\Omega(M,E)$ is locally generated by its
elements of degrees $0$ and $1$, these two derivations are equal.
\end{nproof}

\begin{nlemma}\sl
Under the same assumptions as those of Lemma 4.3.1, let us set,
for each smooth section $X$ of $(E,\tau,M)$,
 \begin{equation*}{\cal L}_\delta(X)
 =\bigl[i(X),\delta\bigr]\,.
 \end{equation*}
Then, for each smooth section $X$ of $(E,\tau,M)$, we have
 \begin{equation*}\bigl[{\cal
L}_\delta(X),\delta\bigr]=\bigl[i(X),\delta^2\bigr]\,;\end{equation*}
for each pair $(X,Y)$ of smooth sections of $(E,\tau,M)$, we have
 \begin{equation*}\bigl[{\cal L}_\delta(X),{\cal L}_\delta(Y)\bigr]
 -{\cal L}_\delta\bigl([X,Y]_\delta\bigr)
 =\Bigl[\bigl[i(X),\delta^2\bigr],i(Y)\Bigr]\end{equation*}
and, for each triple $(X,Y,Z)$ of smooth sections of $(E,\tau,M)$,
we have
 \begin{equation*}i\Bigl(\bigl[X,[Y,Z]_\delta\bigr]_\delta+
          \bigl[Y,[Z,X]_\delta\bigr]_\delta+
          \bigl[Z,[X,Y]_\delta\bigr]_\delta\Bigr)
    =\Bigl[\bigl[[i(X),\delta^2\bigr],
    i(Y)\bigr],i(Z)\Bigr]\,.\end{equation*}
\end{nlemma}

\begin{nproof}
Let us use the graded Jacobi identity. We may write
 \begin{equation*}\bigl[{\cal L}_\delta(X),\delta\bigr]
   =\bigl[[i(X),\delta],\delta\bigr]
   =-\bigl[[\delta,\delta],i(X)\bigr]
    -\bigl[[\delta,i(X)],\delta\bigr]\,.\end{equation*}
Since $[\delta,\delta]=2\delta^2$, we obtain
 \begin{equation*}2\bigl[{\cal
L}_\delta(X),\delta\bigr]=-2\bigl[\delta^2,i(X)\bigr]
    =2\bigl[i(X),\delta^2\bigr]\,,\end{equation*}
which proves the first equality. Similarly, by using again the
graded Jacobi identity and the equality just proven,
  \begin{equation*}
  \begin{split}
  \bigl[{\cal L}_\delta(X),{\cal L}_\delta(Y)\bigr]
  &=\bigl[{\cal L}_\delta(X),[i(Y),\delta]\bigr]\\
  &=-\bigl[i(Y),[\delta,{\cal L}_\delta(X)]\bigr]
      +\bigl[\delta,[{\cal L}_\delta(X),i(Y)]\bigr]\\
  &=-\bigl[[\delta,{\cal L}_\delta(X)],i(Y)\bigr]
      +\bigl[[{\cal L}_\delta(X),i(Y)],\delta\bigr]\\
  &=\bigl[[{\cal L}_\delta(X),\delta],i(Y)\bigr]
      +\bigl[i\bigl([X,Y]_\delta\bigr),\delta\bigr]\\
  &=\bigl[[i(X),\delta^2],i(Y)\bigr]
      +{\cal L}_\delta\bigl([X,Y]_\delta\bigr)\,.
 \end{split}
  \end{equation*}
The second formula is proven. Finally,
  \begin{equation*}
  \begin{split}
  i\Bigl(\bigl[X,[Y,Z]_\delta\bigr]_\delta\Bigr)
  &=\bigl[{\cal L}_\delta(X),i\bigl([Y,Z]_\delta\bigr)\bigr]\\
  &=\bigl[{\cal L}_\delta(X),[{\cal L}_\delta(Y),i(Z)]\bigr]\\
  &=-\bigl[{\cal L}_\delta(Y),[i(Z),{\cal L}_\delta(X)]\bigr]
    -\bigl[i(Z),[{\cal L}_\delta(X),{\cal L}_\delta(Y)]\bigr]\\
  &=\bigl[{\cal L}_\delta(Y),[{\cal L}_\delta(X),i(Z)]\bigr]
    -\bigl[i(Z),{\cal L}_\delta\bigl([X,Y]_\delta\bigr)\bigr]\\
  &\quad-\Bigl[i(Z),\bigl[[i(X),\delta^2],i(Y)\bigr]\Bigr]\\
  &=i\bigl(\bigl[Y,[X,Z]_\delta\bigr]_\delta\bigr)
     +i\bigl(\bigl[[X,Y]_\delta,Z\bigr]_\delta\bigr)\\
  &\quad+\Bigl[\bigl[[i(X),\delta^2],
    i(Y)\bigr],i(Z)\Bigr]\,.
 \end{split}
 \end{equation*}
The proof is complete.
\end{nproof}

\begin{ntheorem}\sl
Let $(E,\tau,M)$ be a vector bundle and $(E^*,\pi,M)$ its dual
bundle. Let $\delta$ be a graded derivation of degree $1$ of the
exterior algebra $\Omega(M,E)$ (notations defined in 3.2.4), which
satisfies
 \begin{equation*}\delta^2=\delta\circ\delta=0\,.\end{equation*}
Then $\delta$ determines a natural Lie algebroid structure on
$(E,\tau,M)$. For that structure, the anchor map $\rho:E\to TM$ is
the unique vector bundle map such that, for each smooth section
$X$ of $(E,\tau,M)$ and each function $f\in C^\infty(M,\RR)$,
 \begin{equation*}i(\rho\circ X)\,df=\langle\delta f,X\rangle\,.\end{equation*}
The bracket $(X,Y)\mapsto\{X,Y\}$ is the $\delta$-bracket defined
in Lemma 4.3.1; it is such that, for each pair $(X,Y)$ of smooth
sections of $(E,\tau,M)$,
 \begin{equation*}i\bigl(\{X,Y\}\bigr)
 =\bigl[[i(X),\delta],i(Y)\bigr]\,.
 \end{equation*}
The $\omega(M,E)$-valued exterior derivative associated to that
Lie algebroid structure (propositions 4.2.1 and 4.2.3) is the
given derivation $\delta$.
\end{ntheorem}

\begin{nproof}
Since $\delta^2=0$, lemmas 4.3.1 and 4.3.2 prove that the
$\delta$-bracket satisfies the Jacobi identity. Let $X$ and $Y$ be
two smooth sections of $(E,\tau,M)$ and $f$ a smooth function on
$M$. By using the definition of the $\delta$-bracket we obtain
 \begin{equation*}i\bigl([X,fY]_\delta\bigr)=f\,i\bigl([X,Y]_\delta\bigr)+
 \bigl({\cal L}_\delta(X)f\bigr)\,i(Y)\,.\end{equation*}
But
 \begin{equation*}{\cal L}_\delta(X)f=\bigl[i(X),\delta\bigr]f=\langle\delta
 f,X\rangle\,,\end{equation*}
since $i(X)f=0$. We must prove now that the value of $\delta(f)$
at any point $x\in M$ depends only on the value of the
differential $df$ of $f$ at that point. We first observe that
$\delta$ being a derivation, the values of $\delta(f)$ in some
open subset $U$ of $M$ depend only on the values of $f$ in that
open subset. Moreover, we have
 \begin{equation*}\delta(1\,f)=\delta(f)=\delta(1)\,f+1\,\delta(f)=\delta(1)\,f+
\delta(f)\,,\end{equation*} which proves that $\delta$ vanishes on
constants.
\par
Let $a\in M$. We use a chart of $M$ whose domain $U$ contains $a$,
and whose local coordinates are denoted by $(x^1,\ldots,x^n)$. In
order to calculate $\delta(f)(a)$, the above remarks allow us to
work in that chart. We may write
 \begin{equation*}f(x)=f(a)+\sum_{i=1}^n(x^i-a^i)\varphi_i(x)\,,\quad\hbox{with}
\quad
        \lim_{x\to a}\varphi_i(x)=\frac{\partial f}{\partial
x^i}(a)\,.\end{equation*} Therefore,
 \begin{equation*}(\delta f)(a)=\sum_{i=1}^n\frac{\partial f}{\partial
x^i}(a)\,\delta(x^i)(a)\,.
 \end{equation*}
We have proven that $\delta(f)(a)$ depends only on $df(a)$, and
that we may write
 \begin{equation*}\delta(f)={}^t\!\rho_\delta\circ df\,,\end{equation*}
where ${}^t\!\rho_\delta:T^*M\to E^*$ is a smooth vector bundle
map. Let $\rho_\delta:E\to TM$ be its transpose. We may now write
  \begin{equation*}[X,fY]_\delta=f[X,Y]_\delta+\langle df,\rho_\delta\circ
X\rangle
  Y\,.\end{equation*}
This proves that the vector bundle $(E,\tau,M)$, with the
$\delta$-bracket and the map $\rho_\delta$ as anchor, is a Lie
algebroid. Finally, by using Propositions 4.2.1 and 4.2.2, we see
that the $\Omega(M,E)$-valued exterior derivative associated to
that Lie algeboid structure is the derivation $\delta$.
\end{nproof}

\subsection{The Schouten-Nijenhuis bracket}

In this subsection $(E,\tau,M,\rho)$ is a Lie algebroid. We have
seen (Propositions 4.1.4 and 4.1.6) that the composition law wich
associates, to each pair $(V,W)$ of sections of the Lie algebroid
$(E,\tau,M,\rho)$, the bracket $\{V,W\}$, extends into a map
$(V,P)\mapsto{\cal L}_\rho(V) P$, defined on $A^1(M,E)\times
A(M,E)$, with values in $A(M,E)$. Theorem~4.4.3 below will show
that this map extends, in a very natural way, into a composition
law $(P,Q)\mapsto[P,Q]$, defined on $A(M,E)\times A(M,E)$, with
values in $A(M,E)$, called the {\it Schouten-Nijenhuis bracket\/}.
That bracket was discovered by Schouten \cite{Schou} for
multivectors on a manifold, and its properties were further
studied by Nijenhuis \cite{Ni}.
\par\smallskip
The following lemmas will be used in the proof of Theorem~4.4.3.

\begin{nlemma}\sl Let $(E,\tau,M,\rho)$ be a Lie algebroid, $p$ and
$q\in\ZZ$, $P\in A^p(M,E)$, $Q\in A^q(M,E)$, $f\in
C^\infty(M,\RR)$ and $\eta\in \Omega(M,E)$. Then
 \begin{equation*}
 \begin{split}
 i(P)\bigl(df\wedge i(Q)\eta\bigr)&-(-1)^pdf\wedge\bigl(i(P)\circ
 i(Q)\eta\bigr)\\
 &+(-1)^{(p-1)q}i(Q)\circ i(P)(df\wedge\eta)\\
 &+(-1)^{(p-1)q+p}i(Q)\bigl(df\wedge i(P)\eta\bigr)\\
 &=0\,.
 \end{split}
 \end{equation*}

\end{nlemma}

\begin{nproof} Let us denote by $E(P,Q,f,\eta)$ the left hand side of
the above equality. We have to prove that $E(P,Q,f,\eta)=0$.
\par
Obviously, $E(P,Q,f,\eta)=0$ when $p<0$, as well as when $q<0$.
When $p=q=0$, we have
 \begin{equation*}E(P,Q,f,\eta)=PQ\,df\wedge\eta-PQ\,df\wedge\eta
 -QP\,df\wedge\eta+QP\,df\wedge\eta=0\,.\end{equation*}
Now we proceed by induction on $p$ and $q$, with the induction
assumption that $E(P,Q,f,\eta)=0$ when $p\leq p_M$ and $q\leq
q_M$, for some integers $p_M$ and $q_M$. Let $P=X\wedge P'$, with
$X\in A^1(M,E)$ and $P'\in A^{p_M}(M,E)$, $Q\in A^q(M,E)$, with
$q\leq q_M$, $f\in C^\infty(M,\RR)$ and $\eta\in\Omega(M,E)$. We
obtain, after some calculations,
 \begin{equation*}
 \begin{split}
 E(P,Q,f,\eta)
 &=E(X\wedge P',Q,f,\eta)\\
 &=(-1)^{p_M+q-1}E\bigl(P',Q,f,i(X)\eta\bigr)\\
 &\quad+(-1)^{p_M}\langle df,X\rangle i(P)\circ i(Q)\eta\\
 &\quad-(-1)^{p_M+p_Mq} \langle df,X\rangle i(Q)\circ i(P)\eta\\
 &=0\,,
 \end{split}
 \end{equation*}
since, by the induction assumption,
$E\bigl(P',Q,f,i(X)\eta\bigr)=0$.
\par
Since every $P\in A^{p_M+1}(M,E)$ is the sum of terms of the form
$X\wedge P'$, with $X\in A^1(M,E)$ and $P'\in A^{p_M}(M,E)$, we
see that $E(P,Q,f,\eta)=0$ for all $p\leq p_M+1$, $q\leq q_M$,
$P\in A^{p_M+1}(M,E)$ and $Q\in A^{q_M}(M,E)$.
\par
Moreover, $P$ and $Q$ play similar parts in $E(P,Q,f,\eta)$, since
we have
 \begin{equation*}E(P,Q,f,\eta)=(-1)^{pq+p+q}E(Q,P,f,\eta)\,.\end{equation*}
Therefore  $E(P,Q,f,\eta)=0$ for all $p\leq p_M+1$, $q\leq q_M+1$,
$P\in A^{p}(M,E)$ and $Q\in A^{q}(M,E)$. By induction we conclude
that  $E(P,Q,f,\eta)=0$ for all $p$ and $q\in\ZZ$, $P\in
A^{p}(M,E)$ and $Q\in A^{q}(M,E)$.
\end{nproof}

\begin{nlemma}\sl Let $(E,\tau,M,\rho)$ be a Lie algebroid, $p$, $q$ and
$r\in\ZZ$, $P\in A^p(M,E)$, $Q\in A^q(M,E)$ and $R\in A^r(M,E)$.
Then
 \begin{equation*}i(R)\circ\bigl[[i(P),d_\rho],i(Q)\bigr]=(-1)^{(p+q-1)r}
 \bigl[[i(P),d_\rho],i(Q)\bigr]\circ i(R)\,.\end{equation*}

\end{nlemma}

\begin{nproof} Let us first assume that $R=V\in A^1(M,E)$. We may write
 \begin{equation*}
 \begin{split}
 i(V)\circ\bigl[[i(P),d_\rho],i(Q)\bigr]
 &=i(V)\circ i(P)\circ d_\rho\circ i(Q)\\
 &\quad-(-1)^pi(V)\circ d_\rho\circ
 i(P)\circ i(Q)\\
 &\quad-(-1)^{(p-1)q}i(V)\circ i(Q)\circ i(P)\circ d_\rho\\
 &\quad+(-1)^{(p-1)q+p} i(V)\circ i(Q)\circ d_\rho\circ i(P)\,.
 \end{split}
 \end{equation*}
We transform the right hand side by pushing the operator $i(V)$
towards the right, using the formulae (proven in 3.2.3~(v) and in
Property~1 of~4.2.3)
 \begin{equation*}i(V)\circ i(P)=(-1)^pi(P)\circ i(V)\quad\hbox{\rm and}\quad
 i(V)\circ d_\rho={\cal L}_\rho(V)-d_\rho\circ i(V)\,.\end{equation*}
We obtain, after rearrangement of the terms,
 \begin{equation*}
 \begin{split}
 i(V)\circ\bigl[[i(P),d_\rho],i(Q)\bigr]
 &=(-1)^{p+q-1} \bigl[[i(P),d_\rho],i(Q)\bigr]\circ i(V)\\
 &\quad+(-1)^pi(P)\circ{\cal L}_\rho(V)\circ i(Q)\\
 &\quad-(-1)^p{\cal L}_\rho(V)\circ i(P)\circ i(Q)\\
 &\quad-(-1)^{(p-1)q+p+q}i(Q)\circ i(P)\circ {\cal L}_\rho(V)\\
 &\quad+(-1)^{(p-1)q+p+q}i(Q)\circ{\cal L}_\rho(V)\circ i(P)\,.
 \end{split}
 \end{equation*}
Now we transform the last four terms of the right hand side by
pushing to the right the operator ${\cal L}_\rho(V)$, using
formulae, proven in 3.2.3~(v) and in Property~4 of~4.1.6, of the
type
 \begin{equation*}i(P)\circ i(Q)=i(P\wedge Q)\quad\hbox{\rm and}\quad
 {\cal L}_\rho(V)\circ i(P)=i(P)\circ{\cal L}_\rho(V)+i\bigl({\cal
 L}_\rho(V)P\bigr)\,.\end{equation*}
The terms containing ${\cal L}_\rho(V)$ become
 \begin{equation*}(-1)^pi\Bigl(P\wedge {\cal L}_\rho(V)Q+\bigl({\cal
 L}_\rho(V)P\bigr)\wedge Q-{\cal L}_\rho(V)(P\wedge Q)\Bigr)\,,\end{equation*}
so they vanish, by Property~3 of~4.1.6. So we have
\begin{equation*}i(V)\circ\bigl[[i(P),d_\rho],i(Q)\bigr]=(-1)^{(p+q-1)}
 \bigl[[i(P),d_\rho],i(Q)\bigr]\circ i(V)\,.\end{equation*}
Now let $R=V_1\wedge \cdots\wedge V_r$ be a decomposable element
in $A^r(M,E)$. Since
 \begin{equation*}i(R)=i(V_1)\circ\cdots\circ i(V_r)\,,\end{equation*}
by using $r$ times the above result, we obtain
\begin{equation*}i(R)\circ\bigl[[i(P),d_\rho],i(Q)\bigr]=(-1)^{(p+q-1)r}
 \bigl[[i(P),d_\rho],i(Q)\bigr]\circ i(R)\,.\end{equation*}
Finally the same result holds for all $R\in A^r(M,E)$ by
linearity. \end{nproof}

\begin{ntheorem}\sl Let $(E,\tau,M,\rho)$ be a Lie algebroid. Let $p$ and
$q\in\ZZ$, and let $P\in A^p(M,E)$, $Q\in A^q(M,E)$. There exists
a unique element in $A^{p+q-1}(M,E)$, called the {\it
Schouten-Nijenhuis bracket\/} of $P$ and $Q$, and denoted by
$[P,Q]$, such that the interior product $i\bigl([P,Q]\bigr)$,
considered as a graded endomorphism of degree $p+q-1$ of the
exterior algebra $\Omega(M,E)$, is given by the formula
 \begin{equation*}i\bigl([P,Q]\bigr)
 =\bigl[[i(P),d_\rho],i(Q)\bigr]\,,
 \end{equation*}
the brackets in the right hand side being the graded brackets of
graded endomorphism (Definition~3.1.3).
\end{ntheorem}

\begin{nproof} We observe that for all $r\in\ZZ$, the map
\begin{equation*}\eta\mapsto
\bigl[[i(P),d_\rho],i(Q)\bigr]\eta\,,\end{equation*} defined on
$\Omega^r(M,E)$, with values in $\Omega^{r-p-q+1}(M,E)$, is
$\RR$-linear. Let us prove that it is in fact
$C^\infty(M,\RR)$-linear. Let $f\in C^\infty(M,\RR)$. By
developing the double graded bracket of endomorphisms, we obtain
after some calculations
 \begin{equation*}
 \begin{split}
 \bigl[[i(P),d_\rho],i(Q)\bigr](f\eta)
 &=f\bigl[[i(P),d_\rho],i(Q)\bigr]\eta\\
 &\quad+i(P)\bigl(df\wedge i(Q)\eta\bigr)-(-1)^pdf\wedge\bigl(i(P)\circ
 i(Q)\eta\bigr)\\
 &\quad+(-1)^{(p-1)q}i(Q)\circ i(P)(df\wedge\eta)\\
 &\quad+(-1)^{(p-1)q+p}i(Q)\bigl(df\wedge i(P)\eta\bigr)\,.
 \end{split}
 \end{equation*}
Lemma~4.4.1 shows that the sum of the last four terms of the right
hand side vanishes, so we obtain
 \begin{equation*}\bigl[[i(P),d_\rho],i(Q)\bigr](f\eta)
 =f\bigl[[i(P),d_\rho],i(Q)\bigr]\eta\,.\end{equation*}
Let us take $r=p+q-1$, and $\eta\in \Omega^{p+q-1}(M,E)$. The map
 \begin{equation*}\eta\mapsto
\bigl[[i(P),d_\rho],i(Q)\bigr]\eta\,,\end{equation*} defined on
$\Omega^{p+q-1}(M,E)$, takes its values in
$\Omega^0(M,E)=C^\infty(M,\RR)$, and is $C^\infty(M,\RR)$-linear.
This proves the existence of a unique element $[P,Q]$ in
$\Omega^{p+q-1}(M,E)$ such that, for all $\eta\in
\Omega^{p+q-1}(M,E)$,
 \begin{equation*}\bigl[[i(P),d_\rho],i(Q)\bigr]\eta
 =i\bigl([P,Q]\bigr)\eta\,.
 \end{equation*}
We still have to prove that the same formula holds for all
$r\in\ZZ$ and all $\eta\in\Omega^r(M,E)$. The formula holds
trivially when $r<p+q-1$, so let us assume that $r>p+q-1$. Let
$\eta\in\Omega^r(M,E)$ and $R\in A^{r-p-q+1}(M,E)$. By using
Lemma~4.4.2, we may write
 \begin{equation*}
 \begin{split}
 i(R)\circ\bigl[[i(P),d_\rho],i(Q)\bigr](\eta)
 &=(-1)^{(p+q-1)(r-p-q+1)}\bigl[[i(P),d_\rho],i(Q)\bigr]\bigl(i(R)\eta\bigr)\\
 &=(-1)^{(p+q-1)(r-p-q+1)}i\bigl([P,Q]\bigr)\bigl(i(R)\eta)\,,
 \end{split}
 \end{equation*}
since $i(R)\eta\in\Omega^{p+q-1}(E)$. Therefore
 \begin{equation*}
 \begin{split}
 i(R)\circ\bigl[[i(P),d_\rho],i(Q)\bigr](\eta)
 &=(-1)^{(p+q-1)(r-p-q+1)}i\bigl([P,Q]\bigr)\circ i(R)\eta)\\
 &=i(R)\circ i\bigl([P,Q]\bigr)\eta\,.
 \end{split}
 \end{equation*}
Since that equality holds for all $\eta\in\Omega^r(M,E)$ and all
$R\in A^{r-p-q+1}(M,E)$, we may conclude that
 \begin{equation*}\bigl[[i(P),d_\rho],i(Q)\bigr]
 =i\bigl([P,Q]\bigr)\,,
 \end{equation*}
and the proof is complete.
\end{nproof}
\par\medskip
In Proposition 4.1.1, we introduced  the Lie derivative with
respect to a section of the Lie algebroid $(E,\tau,M,\rho)$. Now we
define, for all $p\in\ZZ$ and $P\in A^p(M,E)$, the Lie derivative
with respect to $P$. The reader will observe that Property~1 of
Proposition~4.2.3 shows that for $p=1$, the following definition
is in agreement with the definition of the Lie derivative with
respect to an element in $A^1(M,E)$ given in~4.1.1.

\begin{ndefi}\sl
Let $(E,\tau,M,\rho)$ be a Lie algebroid, $p\in\ZZ$ and $P\in
A^p(M,E)$. The {\it Lie derivative\/} with respect to $P$ is the
graded endomorphism of $\Omega(M,P)$, of degree $1-p$, denoted by
${\cal L}_\rho(P)$,
 \begin{equation*}{\cal L}_\rho(P)=\bigl[i(P),d_\rho\bigr]=i(P)\circ d_\rho
 -(-1)^pd_\rho\circ i(P)\,.\end{equation*}

\end{ndefi}

\begin{nrmk}\rm Under the assumptions of Theorem~4.4.3, the above
Definition allows us to write
  \begin{equation*}i\bigl([P,Q]\bigr)=\bigl[{\cal L}_\rho(P),i(Q)\bigr]
  ={\cal L}_\rho(P)\circ i(Q)-(-1)^{(p-1)q}i(Q)\circ{\cal
  L}_\rho(P)\,.\end{equation*}
For $p=1$ and $P=V\in A^1(M,E)$,  this formula is simply Property
4 of Proposition 4.1.6, as shown by the following Proposition.
\end{nrmk}

\begin{nprop}\sl Under the assumptions of Theorem~4.4.3, let
$p=1$, $P=V\in A^1(M,E)$ and $Q\in A^q(M,E)$. The
Schouten-Nijenhuis bracket $[V,Q]$ is simply the Lie derivative of
$Q$ with respect to $V$, as defined in Proposition~4.1.5:
 \begin{equation*}[V,Q]={\cal L}_\rho(V)Q\,.\end{equation*}
\end{nprop}

\begin{nproof} As seen in Remark~4.4.5, we may write
 \begin{equation*}i\bigl([V,Q]\bigr)=\bigl[{\cal L}_\rho(V),i(Q)\bigr]
  ={\cal L}_\rho(V)\circ i(Q)-i(Q)\circ{\cal
  L}_\rho(V)\,.\end{equation*}
Property~4 of Proposition~4.1.6 shows that
 \begin{equation*}i\bigl({\cal L}_\rho(V)Q\bigr)={\cal L}_\rho(V)\circ i(Q)
 -i(Q)\circ{\cal L}_\rho(V)\,.\end{equation*}
Therefore,
 \begin{equation*}i\bigl([V,Q]\bigr)=i\bigl({\cal
L}_\rho(V)Q\bigr)\,,\end{equation*} and finally
 \begin{equation*}[V,Q]={\cal L}_\rho(V)Q\,,\end{equation*}
which ends the proof.
\end{nproof}

\begin{nrmks}\rm\hfill
\par\nobreak\smallskip\noindent
{\rm(i)}\quad{\it The Lie derivative of elements in
$A^p(M,E)$\/.}\quad One may think to extend the range of
application of the Lie derivative with respect to a multivector
$P\in A^p(M,E)$ by setting, for all $q\in\ZZ$ and $Q\in A^q(M,E)$,
 \begin{equation*}{\cal L}_\rho(P)Q=[P,Q]\,,\end{equation*}
the bracket in the right hand side being the Schouten-Nijenhuis
bracket. However, we will avoid the use of that notation because
it may lead to confusions: for $p>1$, $P\in A^p(M,E)$, $q=0$ and
$Q=f\in A^0(M,E)=C^\infty(M,\RR)$, the Schouten-Nijenhuis bracket
$[P,f]$ is an element in $A^{p-1}(M,E)$ which does not vanish in
general. But $f$ can be considered also as an element in
$\Omega^0(M,E)$, and the Lie derivative of $f$ with respect to
$P$, in the sense of Definition~4.4.4, is an element in
$\Omega^{-(p-1)}(M,E)$, therefore vanishes identically. So it
would not be a good idea to write ${\cal L}_\rho(P)f=[P,f]$.
\par\smallskip\noindent
{\rm(ii)}\quad{Lie derivatives and derivations\/.}\quad We have
seen (Property~3 of~4.1.4) that the Lie derivative ${\cal
L}_\rho(V)$ with respect to a section $V\in A^1(M,R)$ of the Lie
algebroid $(E,\tau,M,\rho)$ is a derivation of degree $0$ of the
exterior algebra $\Omega(M,E)$. For $p>1$ and $P\in A^p(M,E)$, the
Lie derivative ${\cal L}_\rho(P)$ is a graded endomorphism of
degree $-(p-1)$ of $\Omega(M,E)$. Therefore, it vanishes
identically on $\Omega^0(M,E)$ and on $\Omega^1(M,E)$. Unless it
vanishes identically, ${\cal L}_\rho(P)$ is not a derivation of
$\Omega(M,E)$.
\end{nrmks}

\begin{nprop}\sl Let $(E,\tau,M,\rho)$ be a Lie algebroid,
$p$ and $q\in\ZZ$, $P\in A^p(M,E)$, $Q\in A^q(M,E)$.
\par\nobreak\smallskip
{\rm 1.} The graded bracket of the Lie derivative ${\cal
L}_\rho(P)$ and the exterior differential $d_\rho$ vanishes
identically:
 \begin{equation*}\bigl[{\cal L}_\rho(P),d_\rho\bigr]={\cal L}_\rho(P)\circ
 d_\rho-(-1)^{p-1}d_\rho\circ{\cal L}_\rho(P)=0\,.\end{equation*}
\par\smallskip
{\rm 2.} The graded bracket of the Lie derivatives ${\cal
L}_\rho(P)$ and ${\cal L}_\rho(Q)$ is equal to the Lie derivative
${\cal L}_\rho\bigl([P,Q]\bigr)$:
 \begin{equation*}\bigl[{\cal L}_\rho(P),{\cal L}_\rho(Q)\bigr]=
 {\cal L}_\rho(P)\circ {\cal L}_\rho(Q)-(-1)^{(p-1)(q-1)}
 {\cal L}_\rho(Q)\circ{\cal L}_\rho(P)={\cal
 L}_\rho\bigl([P,Q]\bigr)\,.\end{equation*}

\end{nprop}

\begin{nproof}
 {\rm 1.} We have seen (3.1.8~(ii)) that the space of graded
endomorphisms of $\Omega(M,E)$, with the graded bracket as
composition law, is a graded Lie algebra. By using the graded
Jacobi identity, we may write
 \begin{equation*}(-1)^p\bigl[[i(P),d_\rho],d_\rho\bigr]
  +(-1)^p\bigl[[d_\rho,d_\rho],i(P)\bigr]
  -\bigl[[d_\rho,i(P)],d_\rho\bigr]=0\,.\end{equation*}
But
 \begin{equation*}[d_\rho,d_\rho]=2 d_\rho\circ d_\rho=0\quad
 \hbox{\rm and}\quad
\bigl[i(P),d_\rho\bigr]=-(-1)^p\bigl[d_\rho,i(P)\bigr]\,.\end{equation*}
So we obtain
 \begin{equation*}2\bigl[[i(P),d_\rho],d_\rho\bigr]=2\bigl[{\cal
 L}_\rho(P),d_\rho\bigr]=0\,.\end{equation*}
\par\smallskip
{\rm 2.} We have
 \begin{equation*}{\cal
 L}_\rho\bigl([P,Q]\bigr)=\bigl[i\bigl([P,Q]\bigr),d_\rho\bigr]
 =\bigl[[{\cal L}_\rho(P),i(Q)],d_\rho\bigr]\,.\end{equation*}
Using the graded Jacobi identity, we may write
 \begin{equation*}
 \begin{split}
 (-1)^{p-1}\bigl[[{\cal L}_\rho(P),i(Q)],d_\rho\bigr]
 &+(-1)^{q(p-1)}\bigl[[i(Q),d_\rho],{\cal L}_\rho(P)\bigr]\\
 &+(-1)^q\bigl[[d_\rho,{\cal L}_\rho(P)],i(Q)\bigr]=0\,.
 \end{split}
 \end{equation*}
But, according to~4.4.4 and Property~1 above,
 \begin{equation*}\bigl[i(Q),d_\rho\bigr]={\cal L}_\rho(Q)\quad\hbox{\rm
 and}\quad\bigl[d_\rho,{\cal L}_\rho(P)\bigr]=0\,.\end{equation*}
So we obtain
 \begin{equation*}{\cal L}_\rho\bigl([P,Q]\bigr)=-(-1)^{(p-1)(q-1)}
 \bigl[{\cal L}_\rho(Q),{\cal L}_\rho(P)\bigr]
 =\bigl[{\cal L}_\rho(P),{\cal L}_\rho(Q)\bigr]\,,\end{equation*}
as announced.
\end{nproof}

\begin{nprop}\sl
 Under the same assumptions as those of
Theorem~4.4.3, the Schouten-Nijenhuis bracket has the following
properties.
 \par\nobreak\smallskip

{\rm 1.} For $f$ and $g\in A^0(M,E)=C^\infty(M,\RR)$,
 \begin{equation*}[f,g]=0\,.\end{equation*}
\par\smallskip

{\rm 2.} For $V\in A^1(M,E)$, $q\in\ZZ$ and $Q\in A^q(M,E)$,
 \begin{equation*}[V,Q]={\cal L}_\rho(V)Q\,.\end{equation*}
\par\smallskip

{\rm 3.} For $V$ and $W \in A^1(M,E)$,
 \begin{equation*}[V,W]=\{V,W\}\,,\end{equation*}
the bracket in the right hand side being the Lie algebroid
bracket.

{\rm 4.} For all $p$ and $q\in\ZZ$, $P\in A^p(M,E)$, $Q\in
A^q(M,E)$,
 \begin{equation*}[P,Q]=-(-1)^{(p-1)(q-1)}[Q,P]\,.\end{equation*}
\par\smallskip

{\rm 5.} Let $p\in\ZZ$, $P\in A^p(M,E)$. The map $Q\mapsto [P,Q]$
is a derivation of degree $p-1$ of the graded exterior algebra
$A(M,E)$. In other words, for $q_1$ and $q_2\in\ZZ$, $Q_1\in
A^{q_1}(M,E)$ and $Q_2\in A^{q_2}(M,E)$,
 \begin{equation*}[P,Q_1\wedge Q_2]=[P,Q_1]\wedge
 Q_2+(-1)^{(p-1)q_1}Q_1\wedge[P,Q_2]\,.\end{equation*}
\par\smallskip

{\rm 6.} Let $p$, $q$ and $r\in\ZZ$, $P\in A^p(M,E)$, $Q\in
A^q(M,E)$ and $R\in A^r(M,E)$. The Schouten-Nijenhuis bracket
satisfies the graded Jacobi identity:
 \begin{equation*}
 \begin{split}
 (-1)^{(p-1)(r-1)}\bigl[[P,Q],R\bigr]
 &+(-1)^{(q-1)(p-1)}\bigl[[Q,R],P\bigr]\\
 &+(-1)^{(r-1)(q-1)}\bigl[[R,P],Q\bigr]\\
 &=0\,.
 \end{split}
 \end{equation*}

\end{nprop}

\begin{nproof}
 {\rm 1.} Let $f$ and $g\in A^0(M,E)$. Then $[f,g]\in
 A^{-1}(M,E)=\{0\}$, therefore $[f,g]=0$.
\par\smallskip

{\rm 2.} See Proposition~4.4.6.
\par\smallskip

{\rm 3.} See Property~1 of Proposition~4.1.6.
\par\smallskip

{\rm 4.}  Let $p$ and $q\in\ZZ$, $P\in A^p(M,E)$, $Q\in A^q(M,E)$.
By using the graded Jacobi identity for graded endomorphisms of
$\Omega(M,E)$, we may write
 \begin{equation*}(-1)^{pq}\bigl[[i(P),d_\rho],i(Q)\bigr]
 +(-1)^{p}\bigl[[d_\rho,i(Q)],i(P)\bigr]
 +(-1)^{q}\bigl[[i(Q),i(P)],d_\rho\bigr]
 =0\,.\end{equation*}
By using
 \begin{equation*}\bigl[i(Q),i(P)\bigr]=i(Q\wedge P)-i(Q\wedge
 P)=0\quad\hbox{and}\quad
 \bigl[d_\rho,i(Q)\bigr]=-(-1)^{q}\bigl[i(Q),d_\rho\bigr]\,,\end{equation*}
we obtain
 \begin{equation*}(-1)^{pq}\bigl[[i(P),d_\rho],i(Q)\bigr]+(-1)^{p+q-1}
 \bigl[[i(Q),d_\rho],i(P)\bigr]=0\,,\end{equation*}
so the result follows immediately.
\par\smallskip

{\rm 5.} Let $p$, $q_1$ and $q_2\in\ZZ$, $P\in A^p(M,E)$, $Q_1\in
A^{q_1}(M,E)$ and $Q_2\in A^{q_2}(M,E)$. We may write
 \begin{equation*}
 \begin{split}
 i\bigl([P,Q_1\wedge Q_2]\bigr)
 &=\bigl[{\cal L}_\rho(P),i(Q_1\wedge Q_2)\bigr]\\
 &={\cal L}_\rho(P)\circ i(Q_1\wedge Q_2)\\
 &\quad-(-1)^{(p-1)(q_1+q_2)}
 i(Q_1\wedge Q_2)\circ{\cal L}_\rho(P)\,.
 \end{split}
 \end{equation*}
We add and substract $(-1)^{(p-1)q_1}i(Q_1)\circ {\cal
L}_\rho(P)\circ i(Q_1)$ from the last expression, and replace
$i(Q_1\wedge Q_2)$ by $i(Q_1)\circ i(Q_2)$. We obtain
 \begin{equation*}i\bigl([P,Q_1\wedge Q_2]\bigr)=\bigl[{\cal
 L}_\rho(P),i(Q_1)\bigr]\circ i(Q_2)
 +(-1)^{(p-1)q_1}i(Q_1)\circ\bigl[{\cal
 L}_\rho(P),i(Q_2)\bigr]\,.\end{equation*}
The result follows immediately.
\par\smallskip

{\rm 6.}
 Let $p$, $q$ and $r\in\ZZ$, $P\in A^p(M,E)$, $Q\in
A^q(M,E)$ and $R\in A^r(M,E)$. By using Property~2 of
Proposition~4.4.8, we may write
 \begin{equation*}
 \begin{split}
 i\bigl(\bigl[[P,Q],R\bigr]\bigr)
 &=\bigl[{\cal L}_\rho([P,Q]),i(R)\bigr]\\
 &=\bigl[[{\cal L}_\rho(P),{\cal L}_\rho(Q)],i(R)\bigr]\,.
 \end{split}
 \end{equation*}
Using the graded Jacobi identity, we obtain
 \begin{equation*}
 \begin{split}
 (-1)^{(p-1)r}\bigl[[{\cal L}_\rho(P),{\cal L}_\rho(Q)],i(R)\bigr]
  &+(-1)^{(q-1)(p-1)}\bigl[[{\cal L}_\rho(Q),i(R)],{\cal
  L}_\rho(P)\bigr]\\
  &+(-1)^{r(q-1)}\bigl[[i(R),{\cal L}_\rho(P)],{\cal L}_\rho(Q)\bigr]
  =0\,.
  \end{split}
  \end{equation*}
But
 \begin{equation*}
 \begin{split}
 \bigl[[{\cal L}_\rho(Q),i(R)],{\cal L}_\rho(P)\bigr]
 &=\bigl[i\bigl([Q,R]\bigr),{\cal L}_\rho(P)\bigr]\\
 &=-(-1)^{(q+r-1)(p-1)}\,\bigl[{\cal
 L}_\rho(P),i\bigl([Q,R]\bigr)\bigr]\\
 &=-(-1)^{(q+r-1)(p-1)}\,i\bigl(\bigl[P,[Q,R]\bigr]\bigr)\\
 &=(-1)^{(q+r-1)(p-1)+(p-1)(q+r-2)}\,i\bigl(\bigl[[Q,R],P\bigr]\bigr)\\
 &=(-1)^{p-1}\,i\bigl(\bigl[[Q,R],P\bigr]\bigr)\,.
 \end{split}
 \end{equation*}
Similarly,
 \begin{equation*}
 \begin{split}
 \bigl[[i(R),{\cal L}_\rho(P)],{\cal L}_\rho(Q)\bigr]
 &=-(-1)^{(p-1)r}\,\bigl[[{\cal L}_\rho(P),i(R)],{\cal
 L}_\rho(Q)\bigr]\\
 &=-(-1)^{(p-1)r}\,\bigl[i\bigl([P,R]\bigr),{\cal
 L}_\rho(Q)\bigr]\\
 &=(-1)^{(p-1)r+(p+r-1)(q-1)}\,\bigl[{\cal
 L}_\rho(Q),i\bigl([P,R]\bigr)\bigr]\\
 &=(-1)^{(p-1)(r+q-1)+r(q-1)}\,i\bigl(\bigl[Q,[P,R]\bigr]\bigr)\\
 &=-(-1)^{(p-1)q+(q-1)(p-2)}\,
 i\bigl(\bigl[[R,P],Q\bigr]\bigr)\\
 &=-(-1)^{p+q}\,i\bigl(\bigl[[R,P],Q\bigr]\bigr)\,.
 \end{split}
 \end{equation*}
Using the above equalities, we obtain
 \begin{equation*}
 \begin{split}
 (-1)^{(p-1)(r-1)}\,i\bigl(\bigl[[P,Q]),R\bigr]\bigr)
 &+(-1)^{(q-1)(p-1)}\,i\bigl(\bigl[[Q,R],P\bigr]\bigr)\\
 &+(-1)^{(r-1)(q-1)}\,i\bigl(\bigl[[R,P],Q\bigr]]\bigr)\\
 &=0\,.
 \end{split}
 \end{equation*}
 The proof is complete.
 \end{nproof}

\begin{nrmks} \rm Let $(E,\tau,M,\rho)$ be a Lie algebroid.
\par\nobreak\smallskip\noindent
{\rm(i)}\quad{\it Degrees for the two algebra structures of
$A(M,E)$\/}.\quad The exterior algebra
$A(M,E)=\oplus_{p\in\ZZ}A^p(M,E)$ of sections of the exterior
powers $(\bigwedge^pE,\tau,M)$, with the exterior product as
composition law, is a graded associative algebra; for that
structure, the space of homogeneous elements of degree $p$ is
$A^p(M,E)$. Proposition~4.4.9 shows that $A(M,E)$, with the
Schouten-Nijenhuis bracket as composition law, is a graded Lie
algebra; for that structure, the space of homogeneous elements of
degree $p$ is not $A^p(M,E)$, but rather $A^{p+1}(M,E)$. For
homogeneous elements in $A(M,E)$, one should therefore make a
distinction between the degree for the graded associative algebra
structure and the degree for the graded Lie algebra structure; an
element in $A^p(M,E)$ has degree $p$ for the graded associative
algebra structure, and degree $p-1$ for the graded Lie algebra
structure.
\par\nobreak\smallskip\noindent
{\rm(ii)}\quad{\it The anchor as a graded Lie algebras
homomorphism\/}.\quad The anchor $\rho:E\to TM$ allows us to
associate to each smooth section $X\in A^1(M,E)$ a smooth vector
field $\rho\circ X$ on $M$; according to Definition 2.1.1, that
correspondence is a Lie algebras homomorphism. We can extend that
map, for all $p\geq 1$, to the space $A^p(M,E)$ of smooth sections
of the $p$-th external power $(\bigwedge^pE,\tau,M)$. First, for a
decomposable element $X_1\wedge\cdots\wedge X_p$, with $X_i\in
A^1(M,E)$, we set
 \begin{equation*}\rho\circ(X_1\wedge\cdots\wedge X_p)=(\rho\circ
 X_1)\wedge\cdots\wedge(\rho\circ X_p)\,.\end{equation*}
For $p=0$, $ f\in A^0(M,E)=C^\infty(M,\RR)$, we set, as a
convention,
  \begin{equation*}\rho\circ f=f\,.\end{equation*}
Then we extend that correspondence to all elements in $A(M,E)$ by
$C^\infty(M,\RR)$-linearity. The map $P\mapsto \rho\circ P$
obtained in that way is a homomorphism from $A(M,E)$ into
$A(M,TM)$, both for their graded associative algebras structures
(with the exterior products as composition laws) and their graded
Lie algebras structures (with the Schouten-Nijenhuis brackets,
associated to the Lie algebroid structure of $(E,\tau,M,\rho)$ and
to the Lie algebroid structure of the tangent bundle
$(TM,\tau_M,M,\id_{TM})$ as composition laws).
\par\smallskip
In 5.2.2~(iii), we will see that when the Lie algebroid under
consideration is the cotangent bundle to a Poisson manifold, the
anchor map has still an additional property: it induces a
cohomology anti-homomorphism.
\end{nrmks}
\

\section{Poisson manifolds and Lie algebroids}
\label{Poisson} In this final section we will show that there
exist very close links between Poisson manifolds and Lie
algebroids.

\subsection{Poisson manifolds}
Poisson manifolds were introduced by A.~Lichnerowicz in the very
important paper \cite{Lich}. Their importance was soon recognozed,
and their properties were investigated in depth by A.~Weinstein
\cite{Wein1}. Let us recall briefly their definition and some of
their properties. The reader is referred to \cite{Lich, Wein1,
Vai} for the proofs of these properties.

\begin{ndefi}\sl
Let $M$ be a smooth manifold. We assume that the space
$C^\infty(M,\RR)$ of smooth functions on $M$ is endowed with a
composition law, denoted by $(f,g)\mapsto\{f,g\}$, for which
$C^\infty(M,\RR)$ is a Lie algebra, which moreover satisfies the
Leibniz-type formula
  \begin{equation*}\{f,gh\}=\{f,g\}h+g\{f,h\}\,.\end{equation*}
We say that the structure defined on $M$ by such a composition law
is a {\it Poisson structure\/}, and that the manifold $M$,
equipped with that structure, is a {\it Poisson manifold\/}.
\end{ndefi}

The following Proposition is due to A.~Lichnerowicz \cite{Lich}.
Independently, A.~Kirillov \cite{Kir} introduced local Lie
algebras (which include both Poisson manifolds and Jacobi
manifolds, which were introduced too by A.~Lichnerowicz
\cite{Lich2}) and obtained, without using the Schouten-Nijenhuis
bracket, an equivalent result and its generalization for Jacobi
manifolds.

\begin{nprop}\sl
On a Poisson manifold $M$, there exists a unique smooth section
$\Lambda\in A^2(M,TM)$, called the {\it Poisson bivector\/}, which
satisfies
 \begin{equation*}[\Lambda,\Lambda]=0\,,\eqno(*)\end{equation*}
such that for any $f$ and $g\in C^\infty(M,\RR)$,
 \begin{equation*}\{f,g\}=\Lambda(df, dg)\,.\eqno(**)\end{equation*}
The bracket in the left hand side of $(*)$ is the
Schouten-Nijenhuis bracket of multivectors on $M$, for the
canonical Lie algebroid structure of $(TM,\tau_M,M)$ (with
$\id_{TM}$ as anchor).
\par
Conversely, let $\Lambda$ be a smooth section of $A^2(TM,M)$. We
use formula $(**)$ to define a composition law on
$C^\infty(M,\RR)$. The structure defined on $M$ by that
composition law is a Poisson structure if and only if $\Lambda$
satisfies formula $(*)$.
\end{nprop}
\par\smallskip
In what follows, we will denote by $(M,\Lambda)$ a manifold $M$
equipped with a Poisson structure whose Poisson bivector is
$\Lambda$.

\subsection{The Lie algebroid structure on the cotangent bundle of
a Poisson manifold}

The next theorem shows that the cotangent bundle of a Poisson
manifold has a canonical structure of Lie algebroid. That property
was discovered by Dazord and Sondaz \cite{DazSon}.

\begin{ntheorem}\sl
Let $(M,\Lambda)$ be a Poisson manifold. There exists, on the
cotangent bundle $(T^*M,\tau_M,M)$, a canonical structure of Lie
algebroid characterized by the following properties:
\par\nobreak\smallskip
\begin{description}
\item{--}\quad the bracket $[\eta,\zeta]$ of two sections $\eta$ and
$\zeta$ of $(T^*M,\tau_M,M)$, {\it i.e.}, of two Pfaff forms on
$M$, is given by the formula
  \begin{equation*}\bigl\langle[\eta,\zeta],X\bigr\rangle
  =\bigl\langle\eta,\bigl[\Lambda,\langle \zeta,X\rangle\bigr]\bigr\rangle
  -\bigl\langle\zeta,\bigl[\Lambda,\langle \eta,X\rangle\bigr]\bigr\rangle
  -[\Lambda,X](\eta,\zeta)\,,\end{equation*}
where $X$ is any smooth vector field on $M$; the bracket in the
right hand side of that formula is the Schouten-Nijenhuis bracket
of multivectors on $M$;

\item{--}\quad the anchor is the vector bundle map
$\Lambda^\sharp:T^*M\to TM$ such that, for each $x\in M$, $\alpha$
and $\beta\in T^*_xM$,

\begin{equation*}\langle\beta,\Lambda^\sharp\alpha\rangle=\Lambda(\alpha,\beta)\
,.\end{equation*}
\end{description}
\end{ntheorem}

\begin{nproof}
We define a linear endomorphism $\delta_\Lambda$ of $A(M,TM)$ by
setting, for each $P\in A(M,TM)$,
 \begin{equation*}\delta_\Lambda(P)=[\Lambda,P]\,,\end{equation*}
where the bracket in the right hand side is the Schouten-Nijenhuis
bracket of multivectors on $M$, {\it i.e.}, the Schouten-Nijenhuis
bracket for the canonical Lie algebroid structure of
$(TM,\tau_M,M)$ (with $\id_{TM}$ as anchor map). When $P$ is in
$A^p(M,TM)$, $\delta_\Lambda(P)$ is in $A^{p+1}(M,TM)$, therefore
$\delta_\Lambda$ is homogeneous of degree $1$. For each $P\in
A^p(M,TM)$ and $Q\in A^q(M,TM)$, we have
 \begin{equation*}
 \begin{split}
 \delta_\Lambda(P\wedge Q)
 &=[\Lambda,P\wedge Q]\\
 &=[\Lambda,P]\wedge Q+(-1)^pP\wedge[\Lambda,Q]\\
 &=\delta_\Lambda(P)\wedge Q+P\wedge\delta_\Lambda(Q)\,.
 \end{split}
 \end{equation*}
This proves that $\delta_\Lambda$ is a graded derivation of degree
$1$ of the exterior algebra $A(M,TM)$.
\par
Moreover, for each $P\in A^p(M,TM)$ we obtain, by using the graded
Jacobi identity,
 \begin{equation*}
 \begin{split}
 \delta_\Lambda\circ \delta_\Lambda(P)
 &=\bigl[\Lambda,[\Lambda,P]\bigr]\\
 &=(-1)^{p-1}\bigl[\Lambda,[P,\Lambda]\bigr]-\bigl[P,[\Lambda,\Lambda]\bigr]\\
 &=-[\Lambda,[\Lambda,P]\bigr]-\bigl[P,[\Lambda,\Lambda]\bigr]\\
 &=-\delta_\Lambda\circ
 \delta_\Lambda(P)-\bigl[P,[\Lambda,\Lambda]\bigr]\,.
 \end{split}
 \end{equation*}
Therefore
 \begin{equation*}2 \delta_\Lambda\circ
 \delta_\Lambda(P)=-\bigl[P,[\Lambda,\Lambda]\bigr]=0\,,\end{equation*}
since $[\Lambda,\Lambda]=0$. We have proven that the graded
derivation $\delta_\Lambda$, of degree $1$, satisfies
 \begin{equation*}\delta_\Lambda^2=\delta_\Lambda\circ
\delta_\Lambda=0\,.\end{equation*} Now we observe that the tangent
bundle $(TM,\tau_M,M)$ can be considered as the dual bundle of the
cotangent bundle $(T^*M,\tau_M,M)$. Therefore, we may apply
Theorem~4.3.3, which shows that there exists on $(T^*M,\tau_M,M)$ a
Lie algebroid structure for which $\delta_M$ is the associated
derivation on the space $\Omega(M,T^*M)=A(M,TM)$ (with the
notations defined in 3.2.4). That theorem also shows that the
bracket of two smooth sections of $(T^*M,\tau_M,M)$, {\it i.e.}, of
two Pfaff forms $\eta$ and $\zeta$ on $M$, is given by the
formula, where $X$ is any smooth vector field on $M$,
 \begin{equation*}\bigl\langle[\eta,\zeta],X\bigr\rangle
 =\bigl\langle\eta,\bigl[\Lambda,\langle \zeta,X\rangle\bigr]\bigr\rangle
 -\bigl\langle\zeta,\bigl[\Lambda,\langle \eta,X\rangle\bigr]\bigr\rangle
 -[\Lambda,X](\eta,\zeta)\,.\end{equation*}
The anchor map $\rho$ is such that, for each $\eta\in
\Omega^1(M,TM)$ and each $f\in C^\infty(M,\RR)$,

\begin{equation*}i(\rho\circ\eta)\,df=\bigl\langle\eta,[\Lambda,f]\bigr\rangle\,
.\end{equation*} The bracket which appears in the right hand sides
of these two formulae is the Schouten-Nijenhuis bracket of
multivectors on $M$. By using Theorem~4.4.3, we see that
 \begin{equation*}[\Lambda,f]=-\Lambda^\sharp(df)\,.\end{equation*}
Therefore,
 \begin{equation*}\langle df,\rho\circ\eta\rangle=i(\rho\circ\eta)\,df
 =\bigl\langle\eta,-\Lambda^\sharp(df)\bigr\rangle=\bigl\langle
 df,\Lambda^\sharp(\eta)\bigr\rangle\,.\end{equation*}
So we have $\rho=\Lambda^\sharp$.
\end{nproof}

\begin{nrmks}\rm Let $(M,\Lambda)$ be a Poisson manifold.
\par\nobreak\smallskip\noindent
{\rm(i)}\quad{\it The bracket of forms of any degrees on
$M$\/}.\quad Since, by Theorem~5.2.1,
$(T^*M,\tau_M,M,\Lambda^\sharp)$ is a Lie algebroid, we can define
a composition law in the space $A(M,T^*M)=\Omega(M,\RR)$ of smooth
differential forms of all degrees on $M$: the Schouten-Nijenhuis
bracket for the Lie algebroid structure of $(T^*M,\tau_M,M)$, with
$\Lambda^\sharp$ as anchor. With that composition law, denoted by
$(\eta,\zeta)\mapsto[\eta,\zeta]$, $\Omega(M,\RR)$ is a graded Lie
algebra. Observe that a form $\eta\in\Omega^p(M,\RR)$, of degree
$p$ for the graded associative algebra structure whose composition
law is the exterior product, has degree $p-1$ for the graded Lie
algebra structure.
\par
The bracket of differential forms on a Poisson manifold was first
discovered for Pfaff forms by Magri and Morosi \cite{MaMo}. It is
related to the Poisson bracket of functions by the formula
 \begin{equation*}[df,dg]=d\{f,g\}\,,\quad\hbox{with}\quad f\ \hbox{and}\ g\in
 C^\infty(M,\RR)\,.\end{equation*}
That bracket was extended to forms of all degrees by Koszul
\cite{Ko2}, and rediscovered, with the Lie algebroid structure of
$T^*M$, by Dazord and Sondaz \cite{DazSon}.
\par\nobreak\smallskip\noindent
{\rm(ii)}\quad{\it The Lichnerowicz-Poisson cohomology\/}.\quad
The derivation $\delta_\Lambda$,
 \begin{equation*}P\mapsto\delta_\Lambda(P)=[\Lambda,P]\,,\quad P\in
A(M,TM)\,,\end{equation*} used in the proof of Theorem 5.2.1, was
first introduced by A.~Lichnerowicz \cite{Lich}, who observed that
it may be used to define a cohomology with elements in $A(M,TM)$
as cochains. He began the study of that cohomology, often called
the Poisson cohomology (but which should be called the
Lichnerowicz-Poisson cohomology). The study of that cohomology was
carried on by Vaisman \cite{Vai}, Huebschmann \cite{Huebs}, Xu
\cite{Xu} and many other authors.
\par\nobreak\smallskip\noindent
{\rm(iii)}\quad{\it The map $\Lambda^\sharp$ as a cohomology
anti-homomorphism\/}.\quad In 4.4.10~(ii), we have seen that the
anchor map $\rho$ of a Lie algebroid $(E,\tau,M,\rho)$ yields a map
$P\mapsto \rho\circ P$ from $A(M,E)$ into $A(M,TM)$, which is both
a homomorphism of graded associative algebras (the composition
laws being the exterior products) and a homomorphism of graded Lie
algebras (the composition laws being the Schouten brackets). When
applied to the Lie algebroid $(T^*M,\tau_M,M,\Lambda^\sharp)$, that
property shows that the map $\eta\mapsto\Lambda^\sharp\circ\eta$
is a homomorphism from the space of differential forms
$\Omega(M,\RR)$ into the space of multivectors $A(M,\RR)$, both
for their structures of graded associative algebras and their
structures of graded Lie algebras. As observed by A.~Lichnerowicz
\cite{Lich}, that map exchanges the exterior derivation $d$ of
differential forms and the derivation $\delta_\Lambda$ of
multivectors (with a sign change, under our sign conventions), in
the following sense: for any $\eta\in \Omega^p(M,\RR)$, we have
 \begin{equation*}\Lambda^\sharp(d\eta)=
 -\delta_\Lambda\bigl(\Lambda^\sharp(\eta)\bigr)
 =-\bigl[\Lambda,\Lambda^\sharp(\eta)\bigr]\,.\end{equation*}
That property is an easy consequence of the formula, valid for any
smooth function $f\in C^\infty(M,\RR)$, which can be derived from
Theorem~4.4.3,
 \begin{equation*}\Lambda^\sharp(df)=-[\Lambda,f]\,.\end{equation*}
The map $\Lambda^\sharp$ therefore induces an anti-homomorphism
from the Lichnerowicz-Poisson cohomology of the Poisson manifold
$(M,\Lambda)$, into its De~Rham cohomology.
\par\nobreak\smallskip\noindent
{\rm(iv)}\quad{\it Lie bialgebroids\/}.\quad Given a Poisson
manifold $(M,\Lambda)$, we have Lie algebroid structures both on
the tangent bundle $(TM,\tau_M,M)$ and on the cotangent bundle
$(T^*M,\tau_M,M)$, with $\id_{TM}:TM\to TM$ and
$\Lambda^\sharp:T^*M\to TM$ as their respective anchor maps.
Moreover, these two Lie algebroid structures are compatible in the
following sense: the derivation $\delta_\Lambda:P\mapsto [\Lambda,
P]$ of the graded associative algebra $A(M,TM)$ (the composition
law being the exterior product) determined by the Lie algebroid
structure of $(T^*M,\tau_M,M)$ is also a derivation for the graded
Lie algebra structure of $A(M,E)$ (the composition law being now
the Schouten-Nijenhuis bracket). We have indeed, as an easy
consequence of the graded Jacobi identity, for $P\in A^p(M,TM)$
and $Q\in A^q(M,TM)$,
 \begin{equation*}
 \begin{split}
 \delta_\Lambda\bigl([P,Q]\bigr)
 &=\bigl[\Lambda,[P,Q]\bigr)
 =\bigl[[\Lambda,P],Q\bigr]+(-1)^{p-1}\bigl[P,[\Lambda,Q]]\\
 &=[\delta_\Lambda P,Q]+(-1)^{p-1}[P,\delta_\Lambda Q]\,.
 \end{split}
 \end{equation*}
When two Lie algebroid structures on two vector bundles in duality
satisfy such a compatibility condition, it is said that that pair
of Lie algebroids is a {\it Lie bialgebroid\/}. The very important
notion of a Lie bialgebroid is dut to K.~Mackenzie and P.~Xu
\cite{MackXu}. Its study was developed by Y.~Kosmann-Schwarzbach
\cite{Kosm} and her student \cite{BanKo} and many other authors.
Recently D.~Iglesias and J.C.~Marrero have introduced a
generalization of that notion in relation with Jacobi manifolds
\cite{IgMar}.
\end{nrmks}

\subsection{The Poisson structure on the dual bundle of a Lie
algebroid} We will now prove that the total space of the dual
bundle $(E^*\pi,M)$ of a Lie algebroid $(E,\tau,M,\rho)$ has a
canonical Poisson structure. This result will allow us to recover
well known results:
\par\smallskip
{--} by taking as Lie algebroid the tangent bundle
$(TM,\tau_M,M,\id_{TM})$ of a smooth manifold $M$, we obtain on
the cotangent bundle $T^*M$ a Poisson structure, which is the
structure associated to its canonical symplectic structure;
\par\smallskip
{--} by taking as Lie algebroid a finite-dimensional Lie algebra
$\cal G$ (the base $M$ being here reduced to a point) we recover,
on its dual ${\cal G}^*$, the canonical Lie-Poisson structure.
\par\smallskip
We will need the following lemma.

\begin{nlemma}\sl
Let $(E,\tau,M,\rho)$ be a Lie algebroid and $(E^*,\pi,M)$ its
dual bundle. To each smooth section $X\in A^1(M,E)$, we associate
a smooth function $\Phi_X$ on $E^*$ by setting
 \begin{equation*}
 \Phi_X(\xi)=\bigl\langle\xi,X\circ\pi(\xi)\bigr\rangle\,.
 \end{equation*}
Let $X\in A^1(M,E)$, $f\in C^\infty(M,\RR)$ and $\xi\in E^*$ be
such that
 \begin{equation*}
 d(\Phi_X+f\circ\pi)(\xi)=0\,.
 \end{equation*}
Then we have
 \begin{equation*}
 X\circ\pi(\xi)=0
 \end{equation*}
and, for any smooth section $Y\in A^1(M,E)$,
 \begin{equation*}
 \Phi_{\{X,Y\}}(\xi)=\langle df,\rho\circ
 Y\rangle\circ\pi(\xi)\,.
 \end{equation*}
\end{nlemma}

\begin{nproof}
Since $d(\Phi_X+f\circ\pi)(\xi)=0$ we have, for any vector
$W\in T_\xi E^*$,
 \begin{equation*}
 \bigl\langle d\Phi_X(\xi),W\bigr\rangle=-\bigl\langle
 df\bigl(\pi(\xi)\bigr),T_\xi\pi(W)\bigr\rangle\,.\eqno(*)
 \end{equation*}
Let us first assume that $W$ is vertical, {\it i.e.}, that
$T_\xi\pi(W)=0$. Then $W$ can be identified with an element of
$\pi^{-1}\bigl(\pi(\xi)\bigr)$, and
$\pi(\xi+sW)=\pi(\xi)$ for all $s\in\RR$. Therefore,
 \begin{equation*}
 \begin{split}
 \bigl\langle d\Phi_X(\xi),W\bigr\rangle
 &=\frac{d}{ds}\Phi_X(\xi+sW)\bigm|_{s=0}\\
 &=\frac{d}{ds}\bigl\langle
   \xi+sW,X\circ\pi(\xi+sW)\bigr\rangle\bigm|_{s=0}\\
 &=\bigl\langle W,X\circ\pi(\xi)\bigr\rangle\\
 &=0\,.
 \end{split}
 \end{equation*}
Since that last equality holds for all
$W\in\pi^{-1}\bigl(\pi(\xi)\bigr)$, we have
$X\circ\pi(\xi)=0$.
\par
Let $(s_1,\ldots,s_k)$ be a family of smooth sections of $\tau$
defined on an open subset $U\subset M$, with $\pi(\xi)\in U$,
such that, for each $y\in U$, $\bigl(s_1(y),\ldots,s_k(y)\bigr)$
is a basis of the fibre $\tau^{-1}(y)$. We have on $U$
 \begin{equation*}
 X=\sum_{i=1}^k X^is_i\,,
 \end{equation*}
where the $X^i$ are smooth functions on $U$ which satisfy
$X_i\circ\pi(\xi)=0$. Then we have, for all $\eta\in
\pi^{-1}(U)$,
 \begin{equation*}
 \Phi_X(\eta)=\sum_{i=1}^k X^i\circ\pi(\eta)
 \bigl\langle\eta,s_i\circ\pi(\eta)\bigr\rangle\,.
 \end{equation*}
Let us now consider a vector $W\in T_\xi E^*$, which may not be
vertical, and a smooth curve $s\mapsto\eta(s)$ in $E^*$ such that
$\eta(0)=\xi$ and $\dfrac{d}{ds}\eta(s)|_{s=0}=W$. We have
 \begin{equation*}
 \begin{split}
 \bigl\langle d\Phi_X(\xi),W\bigr\rangle
 &=\frac{d}{ds}\Phi_X\bigl(\eta(s)\bigr)\bigm|_{s=0}\\
 &=\sum_{i=1}^k\frac{d}{ds}\bigl(X^i\circ\pi\bigl(
   \eta(s)\bigr)\bigr)\bigm|_{s=0}\bigl\langle\xi,
   s_i\circ\pi(\xi)\bigr\rangle\\
 &\quad+
     \sum_{i=1}^k X^i\circ\pi(\xi)\frac{d}{ds}
     \bigl\langle\eta(s),
     s_i\circ\pi\bigl(\eta(s)\bigr)\bigr\rangle
     \bigm|_{s=0}\,.
  \end{split}
 \end{equation*}
Since $X_i\circ\pi(\xi)=0$, the last terms vanish, so we obtain
 \begin{equation*}
 \bigl\langle d\Phi_X(\xi),W\bigr\rangle
 =\sum_{i=1}^k\bigl\langle d(X^i\circ\pi)(\xi),
                           T_\xi\pi(W)
              \bigr\rangle
              \bigl\langle\xi,
                           s_i\circ\pi(\xi)
              \bigr\rangle\,.
 \end{equation*}
Comparing Equation $(*)$ with that last equation, we obtain
 \begin{equation*}
 -\bigl\langle
 df\bigl(\pi(\xi)\bigr),T_\xi\pi(W)\bigr\rangle
 =\sum_{i=1}^k\bigl\langle d(X^i\circ\pi)(\xi),
                           T_\xi\pi(W)
              \bigr\rangle
              \bigl\langle\xi,
                           s_i\circ\pi(\xi)
              \bigr\rangle\,.
 \end{equation*}
Since that equality holds for all $W\in T_\xi E^*$, we have
 \begin{equation*}
 df\circ\pi(\xi)=-\sum_{i=1}^k
 \bigl\langle\xi,s_i\circ\pi(\xi)\bigr\rangle
 d(X^i\circ\pi)(\xi)\,.\eqno(**)
 \end{equation*}
Now, for any smooth section $Y\in A^1(M,E)$, we have on $U$
 \begin{equation*}
 \{X,Y\}=\sum_{i=1}^kX^i\{s_i,Y\}-\sum_{i=1}^k\langle
 dX^i,\rho\circ Y\rangle s_i\,.
 \end{equation*}
Since $X^i\circ\pi(\xi)=0$, we obtain
 \begin{equation*}
 \Phi_{\{X,Y\}}(\xi)=-\sum_{i=1}^k
 \bigl\langle\xi,s_i\circ\pi(\xi)\bigr\rangle
 \bigl\langle d(X^i\circ\pi)(\xi),\rho\circ
 Y\bigl(\pi(\xi)\bigr)\bigr\rangle\,.
 \end{equation*}
From Equation $(**)$ and that last equation, we finally obtain
 \begin{equation*}
 \Phi_{\{X,Y\}}(\xi)=\langle df,\rho\circ
 Y\rangle\circ\pi(\xi)\,
 \end{equation*}
and our proof is complete.
\end{nproof}

\begin{ntheorem}\sl
Let $(E,\tau,M,\rho)$ be a Lie algebroid and $(E^*,\pi,M)$ its
dual bundle. To each smooth section $X\in A^1(M,E)$, we associate
a smooth function $\Phi_X$ on $E^*$ by setting
 \begin{equation*}
 \Phi_X(\xi)=\bigl\langle\xi,X\circ\pi(\xi)\bigr\rangle\,.
 \end{equation*}
There exists on $E^*$ a unique Poisson structure such that, for
each pair $(X,Y)$ of smooth sections of $\tau$,
 \begin{equation*}
 \{\Phi_X,\Phi_Y\}=\Phi_{\{X,Y\}}\,,
 \end{equation*}
the bracket in the left hand side being the Poisson bracket of
functions on $E^*$, and the bracket in the right hand side the
bracket of sections of the Lie algebroid $(E,\tau,M,\rho)$.
\end{ntheorem}

\begin{nproof}
If such a Poisson structure exists, it must be such that for any
$X\in A^1(M,E)$, $f$ and $g\in C^\infty(M,\RR)$,
 \begin{equation*}
 \{\Phi_X,g\circ\pi\}=\bigl({\cal L}(\rho\circ
 X)g\bigr)\circ\pi\,,\qquad
 \{f\circ\pi,g\circ\pi\}=0\,.
 \end{equation*}
These identities are consequences of the property of the Lie
algebroid bracket:
 \begin{equation*}
 \{fX,gY\}=fg\{X,Y\}+\bigl(f{\cal L}(\rho\circ X)g\bigr)Y
                    -\bigl(g{\cal L}(\rho\circ Y)f\bigr)X\,,
 \end{equation*}
which implies
 \begin{equation*}
 \Phi_{\{fX,gY\}}=(fg\circ\pi)\Phi_{\{X,Y\}}
                  +\bigl(f{\cal L}(\rho\circ X)g\bigr)\circ\pi\,\Phi_Y
                  -\bigl(g{\cal L}(\rho\circ
                  Y)f\bigr)\circ\pi\,\Phi_X\,.
 \end{equation*}
Now we observe that for any $\xi\in E^*$, $\eta_1$ and $\eta_2\in
T^*_\xi E^*$, there exist (non unique) pairs $(X_1,f_1)$ and
$(X_2,f_2)$, with $X_1$ and $X_2\in A^1(M,E)$, $f_1$ and $f_2\in
C^\infty(M,\RR)$, such that
 \begin{equation*}
 \eta_1=d(\Phi_{X_1}+f_1\circ\pi)(\xi)\,,\quad
 \eta_2=d(\Phi_{X_2}+f_2\circ\pi)(\xi)\,.
 \end{equation*}
The Poisson bivector $\Lambda$ on $E^*$ is therefore given by the
equation
 \begin{equation*}
 \begin{split}
 \Lambda(\xi)(\eta_1,\eta_2)
 &=\{\Phi_{X_1}+f_1\circ\pi,
 \Phi_{X_2}+f_2\circ\pi\}(\xi)\\
 &=\Phi_{\{X_1,X_2\}}(\xi)+\bigl({\cal L}(\rho\circ
 X_1)f_2\bigr)\circ\pi(\xi)-
 \bigl({\cal L}(\rho\circ
 X_2)f_1\bigr)\circ\pi(\xi)\,.
 \end{split}
 \end{equation*}
This proves that if such a Poisson structure exists, it is unique.
In order to prove its existence we still have to check that the
right hand side of the above formula depends only on $\eta_1$ and
$\eta_2$, not on the particular choices we have made for
$(X_1,f_1)$ and $(X_2,f_2)$. Since that member depends bilinearly
on $(X_1,f_1)$ and $(X_2,f_2)$ and is skew-symmetric, it is enough
to prove that if $X_1\in A^1(M,E)$ and $f_1\in C^\infty(M,\RR)$
are such that
 \begin{equation*}
 d(\Phi_{X_1}+f_1\circ\pi)(\xi)=0\,,
 \end{equation*}
then for any $X_2\in A^1(M,E)$ and $f_2\in C^\infty(M,\RR)$,
 \begin{equation*}
 \Phi_{\{X_1,X_2\}}(\xi)+\bigl({\cal L}(\rho\circ
 X_1)f_2\bigr)\circ\pi(\xi)-
 \bigl({\cal L}(\rho\circ
 X_2)f_1\bigr)\circ\pi(\xi)=0\,.
 \end{equation*}
But that equality follows immediately from Lemma~5.3.1. We have
proven the existence of the bivector $\Lambda$, and the formula
used for its definition proves that it is smooth.
\par
Finally, we observe that the Poisson bracket defined by means of
$\Lambda$, when restricted to functions on $E^*$ of the type
$\Phi_X+f\circ\pi$, with $X\in A^1(M,E)$ and $f\in
C^\infty(M,\RR)$, satisfies the Jacobi identity. Therefore, for
any $\xi\in E^*$, $\eta_1$, $\eta_2$ and $\eta_3\in T^*_\xi E^*$
which are the differentials, at $\xi$, of functions of the type
$\Phi_X+f\circ\pi$, the Schouten bracket $[\Lambda,\Lambda]$
satisfies $[\Lambda,\Lambda](\xi)(\eta_1,\eta_2,\eta_3)=0$. Since
all elements in $T^*_\xi E^*$ are differentials, at $\xi$, of
functions of the type $\Phi_X+f\circ\pi$, we have proven the
$[\Lambda,\Lambda]$ vanishes identically, in other words, that
$\Lambda$ is a Poisson bivector.
\end{nproof}

\begin{nprop}\sl
Let $(E,\tau,M,\rho)$ be a Lie algebroid and $(E^*,\pi,M)$ its
dual bundle. The Poisson structure on $E^*$ defined in
Theorem~5.3.2 has the following properties:
\par\smallskip\noindent
{\rm 1.}\quad For any $X\in A^1(M,E)$, $f$ and $g\in
C^\infty(M,\RR)$,
 \begin{equation*}
 \{\Phi_X,g\circ\pi\}=\bigl({\cal L}(\rho\circ
 X)g\bigr)\circ\pi\,,\qquad
 \{f\circ\pi,g\circ\pi\}=0\,,
 \end{equation*}
where $\Phi_X$ is the smooth function on $M$ associated to the
smooth section $X$ as indicated in Theorem~5.3.2.
\par\smallskip\noindent
{\rm 2.}\quad The transpose $\,{}^t\!\rho:T^*M\to E^*$ of the
anchor map $\rho:E\to TM$ is a Poisson map (the cotangent bundle
being endowed with the Poisson structure associated to its
canonical symplectic structure).
\par\smallskip\noindent
{\rm 3.}\quad The Poisson structure on $E^*$ is homogeneous with
respect to the vector field $Z$ which generates the homotheties in
the fibres of $(E^*,\pi,M)$, which means that its Poisson
bivector $\Lambda$ satisfies
 \begin{equation*}
 [Z,\Lambda]=-\Lambda\,.
 \end{equation*}
\end{nprop}

\begin{nproof}
We have proven Properties 1 in the proof of Theorem~5.3.2. In
order to prove Property 2, we must prove that for all pairs
$(h_1,h_2)$ of smooth functions on $E^*$,
 \begin{equation*}
 \{h_1\circ\,{}^t\!\rho,h_2\circ\,{}^t\!\rho\}
 =\{h_1,h_2\}\circ\,{}^t\!\rho\,,
 \end{equation*}
the bracket in the left hand side being the Poisson bracket of
functions on $T^*M$, and the bracket in the right hand side the
Poisson bracket of functions on $E^*$. It is enough to check that
property when $h_1$ and $h_2$ are of the type $\Phi_X$, where
$X\in A^1(M,E)$, or of the type $f\circ\pi$, with $f\in
C^\infty(M,\RR)$, since the differentials of functions of these
two types generate $T^*E^*$. For $h_1=\Phi_X$ and $h_2=\Phi_Y$,
with $X$ and $Y\in A^1(M,E)$, and $\zeta\in T^*M$, we have
 \begin{equation*}
 \begin{split}
 \{\Phi_X,\Phi_Y\}\circ\,{}^t\!\rho(\zeta)
 &=\Phi_{\{X,Y\}}\circ\,{}^t\!\rho(\zeta)\\
 &=\bigl\langle\,{}^t\!\rho(\zeta),\{X,Y\}
   \circ\pi\circ\,{}^t\!\rho(\zeta)\bigr\rangle\\
 &=\bigl\langle\zeta,\rho\circ\{X,Y\}\circ\tau_M(\zeta\bigr\rangle\\
 &=\bigl\langle\zeta,[\rho\circ X,\rho\circ Y]\circ\tau_M(\zeta\bigr\rangle\,,
 \end{split}
 \end{equation*}
since the canonical projection $\tau_M:T^*M\to M$ satisfies
$\pi\circ\,{}^t\!\rho=\tau_M$. But let us recall a well known
property of the Poisson bracket of functions on $T^*M$
(\cite{LibMa}, exercise 17.5 page 182). To any vector field
$\widehat X$ on $M$, we associate the function $\Psi_{\widehat X}$
on $T^*M$ by setting, for each $\zeta\in T^*M$,
 \begin{equation*}
 \Psi_{\widehat X}(\zeta)=\bigl\langle\zeta,\widehat
 X\circ\tau_M(\zeta)\bigr\rangle\,.
 \end{equation*}
Then, for any pair $(\widehat X,\widehat Y)$ of vector fields on
$M$,
 \begin{equation*}
 \{\Psi_{\widehat X},\Psi_{\widehat Y}\}=\Psi_{[\widehat X,\widehat Y]}\,.
 \end{equation*}
By using $\tau_M=\pi\circ\,{}^t\!\rho$, we easily see that for
each $X\in A^1(M,E)$,
 \begin{equation*}
 \Psi_{\rho\circ X}=\Phi_X\circ\,{}^t\!\rho\,.
 \end{equation*}
Returning to our pair of sections $X$ and $Y\in A^1(M,E)$, we see
that
 \begin{equation*}
 \{\Phi_X,\Phi_Y\}\circ\,{}^t\!\rho(\zeta)=\Psi_{[\rho\circ
 X,\rho\circ Y]}(\zeta)=\{\Psi_{\rho\circ X},\Psi_{\rho\circ
 Y}\}(\zeta)
 =\{\Phi_X\circ\,{}^t\!\rho,\Phi_Y\circ\,{}^t\!\rho\}(\zeta)\,.
 \end{equation*}
Now for $h_1=\Phi_X$ and $h_2=f\circ\pi$, with $X\in A^1(M,E)$
and $f\in C^\infty(M,\RR)$, we have
 \begin{equation*}
 \begin{split}
 \{\Phi_X\circ\,{}^t\!\rho,f\circ\pi\circ\,{}^t\!\rho\}
 &=\{\Psi_{\rho\circ X},f\circ\tau_M\}\\
 &={\cal L}(\rho\circ X)f\circ\tau_M\\
 &={\cal L}(\rho\circ X)f\circ\pi\circ\,{}^t\!\rho\\
 &=\{\Phi_X,f\circ\pi\}\circ\,{}^t\!\rho\,.
 \end{split}
 \end{equation*}
Similarly, for $h_1=f\circ\pi$ and $h_2=g\circ\pi$, we have
 \begin{equation*}
 \{f\circ\pi,g\circ\pi\}\circ\,{}^t\!\rho=0=\{f\circ\tau_M,g\circ\tau_M\}
 =\{f\circ\pi\circ\,{}^t\!\rho,g\circ\pi\circ\,{}^t\!\rho\}\,.
 \end{equation*}
Property 2 is proven. Let us now prove Property~3. The vector
field $Z$ on the vector bundle $(E^*,\pi,M)$ generates the
homotheties in the fibres, therefore its reduced flow is
$(t,\xi)\mapsto H_t(\xi)=e^t\xi$, with $t\in\RR$, $\xi\in E^*$.
For any smooth section $X\in A^1(M,E)$ and any $t\in\RR$, we have
 \begin{equation*}
 (H_t^*\Phi_X)(\xi)=\Phi_X\circ H_t(\xi)=e^t\Phi_X(\xi)\,.
 \end{equation*}
Therefore, for $X$ and $Y\in A^1(M,E)$,
 \begin{equation*}
 H_t^*\bigl(\Lambda(d\Phi_X,d\Phi_Y)\bigr)
 =\{\Phi_X,\Phi_Y\}\circ H_t=\Phi_{\{X,Y\}}\circ H_t
 =e^t\Phi_{\{X,Y\}}\,.
 \end{equation*}
But we may also write
 \begin{equation*}
 H_t^*\bigl(\Lambda(d\Phi_X,d\Phi_Y)\bigr)
 =(H_t^*\Lambda)\bigl(H_t^*d\Phi_X,H_t^*d\Phi_Y\bigr)
 =e^{2t}(H_t^*\Lambda)\bigl(d\Phi_X,d\Phi_Y\bigr)
 \end{equation*}
Since the differentials of functions of the type $\Phi_X$ generate
$T^*E$, except along its zero section, this result proves that,
except maybe along the zero section, we have
 \begin{equation*}
 H_t^*\bigl(\Lambda)=e^{-t}\Lambda\,.
 \end{equation*}
By continuity, hat equality holds everywhere. Finally,
 \begin{equation*}
 [Z,\Lambda]={\cal
 L}(Z)\Lambda=\frac{d}{dt}H_t^*(\Lambda)\bigm|_{t=0}=-\Lambda\,,
 \end{equation*}
and our proof is complete.
\end{nproof}

\

%% ----------------------- REFERENCES -------------------------------%%


\begin{thebibliography}{99}

\bibitem{Albert} C.~Albert, P.~Dazord, {\it Th\'{e}orie des groupo\"{\i}des
symplectiques, Ch.~II\/}, Publ. D\'{e}pt. Math., Univ. Lyon I,
nouvelle s\'{e}rie (1990), 27--99.

\bibitem{BanKo} M.~Bangoura and Y.~Kosmann-Schwarzbach, {\it
\'Equation de Yang-Baxter dynamique classique et alg\'{e}bro\"{\i}de de
Lie\/}, C.R. Acad. Sci. Paris 327, I, 1998, 541--546.

\bibitem{CannW} A.~Cannas da Silva and A.~Weinstein, {\it
Geometric models for noncommutative algebras\/}, Berkeley
Mathematics Lecture Notes 10, Amer. Math. Soc., Providence, 1999.

\bibitem{Coste} A.~Coste, P.~Dazord and A.~Weinstein, {\it
Groupo\"{\i}des symplectiques\/}, Publ. D\'{e}pt. Math., Univ. Lyon I, 2/A
(1987), 1--64.

\bibitem{DazSon} P.~Dazord and D.~Sondaz, {\it Vari\'{e}t\'{e}s de
Poisson, Alg\'{e}bro\"{\i}des de Lie\/}, Publ. D\'{e}pt. Math., Univ. Lyon I,
nouvelle s\'{e}rie 1/B (1988), 1--68.

\bibitem{GrHV} W.~Greub, S.~Halperin and R.~Vanstone, {\it
Connections, curvature and cohomology\/}, vol I, Academic Press,
New York, London, 1972.

\bibitem{Huebs} J.~Huebschmann, {\it Poisson cohomology and
quantization\/}, J. Reine Angew. Math.~408 (1990), 57--113.

\bibitem{IgMar} D.~Iglesias and J.C.~Marrero, {\it Generalized Lie
bialgebroids and Jacobi structures\/}, preprint, Universidad de la
Laguna, Spain, 2000.

\bibitem{Kara} M.~Karasev, {\it Analogues of the objects of Lie
group theory for nonlinear Poisson brackets\/}, Math.~USSR
Izvest.~28 (1987), 497--527.

\bibitem{Kir} A.~Kirillov, {\it Local Lie algebras\/}, Russian Math.~Surveys
31 (1976), 55-75.

\bibitem{Kosm} Y.~Kosmann-Schwarzbach, {\it Exact Gerstenhaber algebras and
Lie bialgebroids\/}, Acta Appl.~Math.~41 (1995), 153-165.

\bibitem{Ko1} J.-L.~Koszul, {\it Homologie et cohomologie des
alg\`{e}bres de Lie\/}, Bull. Soc. Math. France 78 (1950), 65--127.

\bibitem{Ko2} J.-L.~Koszul, {\it Crochet de Schouten-Nijenhuis et
cohomologie\/}, Ast\'{e}risque, hors s\'{e}rie, 1985, 257--271.

\bibitem{Lib1} P.~Libermann, {\it On symplectic and contact
groupoids\/}, Differential Geometry and its applications, Proc.
Conf. Opava, August 24--28, 1992, Silesian University, Opava,
1993, 29--45.

\bibitem{Lib2} P.~Libermann, {\it Lie algebroids and mechanics\/},
Archivum mathematicum 32 (1996), 1147--162.

\bibitem{LibMa} P.~Libermann and Ch.-M.~Marle, {\it Symplectic Geometry and
Analytical Mechanics\/}, Kluwer, Dordrecht, 1987.

\bibitem{Lich} A.~Lichnerowicz, {\it Les vari\'{e}t\'{e}s de Poisson et leurs
alg\`{e}bres de Lie associ\'{e}es\/}, J. Differential Geometry 12 (1977),
253-300.

\bibitem{Lich2} A.~Lichnerowicz, {\it Les vari\'{e}t\'{e}s de Jacobi  et leurs
alg\`{e}bres de Lie associ\'{e}es\/}, J. Math. pures et appl.~57 (1978),
453--488.

\bibitem{Mack} K.C.H.~Mackenzie, {\it Lie groupoids and Lie algebroids
in differential geometry\/}, London Math.~Soc.~Lecture notes
series 124, Cambridge University Press, Cambridge (1987).

\bibitem{MackXu} K.C.H.~Mackenzie and P.~Xu, {\it Lie bialgebroids and Poisson
groupoids\/}, Duke Math. J. 73 (1994), 415--452.

\bibitem{MaMo} F.~Magri and C.~Morosi, {\it A geometrical
characterization of integrable Hamiltonian systems through the
theory of Poisson-Nijenhuis manifolds\/}, Quaderno S. 19 (1984),
Universit\`a di Milano.

\bibitem{Ni} A.~Nijenhuis, {\it Jacobi-type identities for
bilinear differential concomitants of certain tensor fields\/},
Indag. Math. 17 (1955), 390--403.

\bibitem{Pra} J.~Pradines, Th\'eorie de Lie pour les groupo\"{\i}des
diff\'erentiables. Calcul diff\'erentiel dans la cat\'egorie des
groupo\"{\i}des infinit\'esimaux, C.R. Acad. Sci. Paris 264 A (1967),
245-248.

\bibitem{Schou} J.~A.~Schouten, {\it On the differential operators of
first order in tensor calculus\/}, Convegno Intern. Geom. Diff.
Italia, Cremonese, Roma (1953), 1--7.

\bibitem{Vai} I.~Vaisman, {\it Lectures on the Geometry of Poisson
Manifolds\/}, Birkh\"{a}user, Basel, Boston, Berlin (1994).

\bibitem{Wein1} A.~Weinstein, {\it The local structure of Poisson
manifolds\/}, J.~Differential Geometry 18 (1983), 523--557.

\bibitem{Wein2} A.~Weinstein, {\it Symplectic groupoids and Poisson
manifolds\/}, Bull.~Amer. Math.~Soc.~16 (1987), 101--103.

\bibitem{Wein3} A.~Weinstein, {\it Groupoids: unifying internal and
external symmetry, a tour through some examples\/}, Notices of the
Amer. Math. Soc.~43 (1996), 744--752.

\bibitem{Xu} P.~Xu, {\it Poisson cohomology of regular Poisson
manifolds\/}, Ann. Inst. Fourier, Grenoble, 42, 4 (1992),
967--988.

\end{thebibliography}
\end{document}